\documentclass[12pt]{amsart}
\usepackage{cond_bootstrap}

\title{}
\date{\today}

\begin{document}

\begin{center}
	\Large \textsc{On the Conditional Consistency of the Bootstrap}
	\bigskip
\end{center}

\vspace*{5ex minus 1ex}
\begin{center}
	Denis Kojevnikov\,\orcidlink{0009-0008-2514-5835} \\
	\medskip
	\textit{Tilburg University}
	\medskip
\end{center}

\begin{abstract}
{\footnotesize This paper focuses on the conditional consistency of the bootstrap. The target is the conditional law of a statistic given a common shock, and the bootstrap approximates it by the law of a resampled statistic conditional on the shock and the data. We call the bootstrap consistent if these laws are close in the Kolmogorov distance for almost every realization of the shock. The conditioning $\sigma$-algebras are allowed to be unrelated across sample sizes, and both the comparison law and the dimension of the data may vary with the sample size. We give sufficient conditions for consistency in this sense, show that they yield asymptotically correct conditional coverage for bootstrap quantiles, and provide a finite-sample bound on the coverage error. We also relate the Kolmogorov formulation to the merging of probability measures. The results are specialized to statistics depending smoothly on an estimator or on an estimated nuisance parameter.}

\bigskip

{\footnotesize \noindent \textsc{Keywords.} Bootstrap consistency; Common shocks; Conditional coverage; Merging of probability kernels; Smooth function model.}
\end{abstract}

\maketitle

\section{Introduction}
\label{sec:intro}

In many settings, the data are driven by a common shock \citep[e.g.,][]{Andrews:05,Kuersteiner/Prucha:13}, which induces dependence across observations that does not vanish as the sample grows. In this case, limit theorems are typically available given the shock, and an unconditional limit is a mixture over the shock, provided that such a limit exists. Since only one realization of the shock is observed, this mixture cannot be recovered from the data, and a resampling scheme, which operates on the realized sample, can approximate only the conditional law.

An example is the class of network dependent processes, for which a law of large numbers and a central limit theorem are available in \cite{Kojevnikov/Marmer/Song:21}. When a single network is observed, inference on the parameters of such a process is naturally conditional on the network or on the environment that generates it. A typical reason to use the bootstrap in this setting is to avoid estimating the long-run variance, which is unreliable when the sample is small relative to the dimension. Without self-normalization, the conditional law of the statistic need not converge. This paper studies the conditional consistency of the bootstrap \citep{Efron:79} without reference to the underlying data-generating mechanism.

We adapt the notions of convergence of random probability measures developed in \cite{Belyaev:00} and \cite{Berti:06}. For bootstrap replicates that are conditionally i.i.d.\ given the data, \citet{Bucher/Kojadinovic:19} show that bootstrap consistency is equivalent to the joint unconditional weak convergence of the statistic and two of its replicates to independent copies of the same limit. They extend the equivalence to processes in $\ell^{\infty}(\csdot)$, where conditional laws need not exist, and show that the limit can be replaced with a fixed probability measure. The convergence notion behind this equivalence is the conditional weak convergence of \citet[Section 1.13]{vanderVaart/Wellner:23}, where a sequence converges in distribution conditionally on sub-$\sigma$-algebras if the bounded-Lipschitz distance between its conditional law and a fixed nonrandom law tends to zero in outer probability or almost uniformly. The equivalence above is the bootstrap interpretation of ``Hoeffding's trick'' (Lemma 1.13.4 therein), whose proof relies on the nonrandomness of the limit law.

The framework introduced below does not require a common limit. We compare the bootstrap and sample distributions at each sample size, in the manner of the \emph{merging} of probability measures \citep{DAristotile/Diaconis/Freedman:88,Davydov/Rotar:09}, allowing, as in \cite{vanderVaart/Wellner:23}, the conditioning $\sigma$-algebras to be unrelated across sample sizes, and allowing the comparison laws to be random. The spaces in which the data and parameters take values can also vary as the sample size grows. Consistency is formulated as the closeness of two conditional laws in the Kolmogorov distance. This distance is taken as the primitive notion because it directly yields the uniform validity of bootstrap quantiles and remains attainable when the comparison sequence varies with $n$.

The paper is organized as follows. Section \ref{sec:cond_conv} provides the theory of convergence for conditional laws against a possibly random comparison that can vary with the sample size. It collects counterparts of the continuous mapping theorem, the Cram\'er--Wold device, and L\'evy's continuity theorem for merging sequences of conditional laws and builds the P\'olya-type bridge between merging and the Kolmogorov distance in both directions. Section \ref{sec:cond_bootstrap} applies this to the bootstrap, providing a bound on the conditional coverage error for each sample size and sufficient conditions for conditional consistency under different modes of convergence. Section \ref{sec:sfm} specializes these results to statistics that depend smoothly on an estimator or on an estimated nuisance parameter. Finally, Section \ref{sec:examples} illustrates the results with maxima in fixed dimension, whose conditional laws cannot be reduced to a fixed comparison by rescaling, and with maxima and correlations in growing dimension.

\subsection*{Preliminaries}
Throughout, $\N$, $\Q$, and $\R$ denote the sets of natural, rational, and real numbers, respectively, $\norm{\csdot}$ is the Euclidean norm, and $\B(\csdot)$ is the Borel $\sigma$-algebra of the indicated space. We write $\X$, $\Y$, and $\Z$ for Polish spaces, choosing complete compatible metrics $\varrho_{\X}$ and $\varrho_{\Y}$ when necessary, and $d\ge 1$ for a generic finite dimension. The collection of compact subsets of $\X$ is denoted by $\K_{\X}$, and similarly for $\Y$ and $\Z$. A sequence $\seq{\mu_n}$ of Borel probability measures on $\X$ \emph{converges weakly} to $\mu$, written as $\mu_n\Rightarrow\mu$, if $\int f\,d\mu_n\to\int f\,d\mu$ for every $f\in \Cb{\X}$, the set of real-valued bounded continuous functions on $\X$.

We take $(\Omega,\H,\PM)$ to be the underlying probability space, for which $\E$ denotes expectation with respect to $\PM$. Given a sequence $\{\F_n\}$ of sub-$\sigma$-algebras of $\H$, we define $\F^{\wedge}\eqdef\bigcap_n\F_n$, $\F_n^{\vee}\eqdef\bigvee_{m\ge n}\F_m$, $\F^{\vee}\eqdef\F_1^{\vee}$, and $\F_{\infty}\eqdef\bigcap_{n\ge1}\F_n^{\vee}$. Also, $\H_{\varnothing}\eqdef\{\varnothing,\Omega\}$, and for any $A$, $B$, and $C$ that are either random elements or $\sigma$-algebras, we write $\CI{A}{B}{C}$ if $\sigma(A)$ and $\sigma(B)$ are conditionally independent given $\sigma(C)$.

A \emph{probability kernel} on $\X$ is a map $\QM:\Omega\times\B(\X)\to[0,1]$ such that (a) $B\mapsto \QM(\omega,B)$ is a probability measure on $(\X,\B(\X))$ for every $\omega\in\Omega$; and (b) $\omega\mapsto \QM(\omega,B)$ is $\H$\-/measurable for every $B\in\B(\X)$. It is \emph{$\F$\-/measurable}, for a sub-$\sigma$-algebra $\F\subset\H$, if $\omega\mapsto\QM(\omega,B)$ is $\F$\-/measurable for every $B$.\footnote{
	Equivalently, $\QM$ is $\F$\-/measurable if and only if $\omega\mapsto\QM(\omega,\csdot)$ is measurable from $(\Omega,\F)$ into the set of probability measures on $\X$ equipped with the Borel $\sigma$-algebra of the weak topology.
}
A \emph{regular conditional distribution} of an $\X$\-/valued random element $X$ given $\F$ is an $\F$\-/measurable probability kernel $\QM_X^{\F}$ on $\X$ such that for every $B\in\B(\X)$, $\QM_X^{\F}(\csdot,B)$ is a version of $\PR{X\in B\mid\F}$. With $\X=\R$, its \textit{regular conditional cdf} $\FD{X}{\F}$ is the cdf of the corresponding kernel, i.e., $\FD{X}{\F}(\omega,x)=\QM_X^{\F}(\omega,(-\infty,x])$.

Let $S\subset\X$. A function $f:\Omega\times S\to \Y$ is called a \emph{Carath\'eodory function} if (a) for each $x\in S$, the map $\omega\mapsto f(\omega,x)$ is $(\H,\B(\Y))$\-/measurable; and (b) for each $\omega\in\Omega$, the map $x\mapsto f(\omega,x)$ is continuous. Every such function is jointly measurable \citep[see, e.g.,][Lemma 4.51]{Aliprantis/Border:06:InfiniteDimensionalAnalysis}. For a probability kernel $\QM$ on $\X$ and a jointly measurable function $g:\Omega\times\X\to\Y$, we write $\QM\circ g^{-1}$ for the probability kernel $\QM'$ on $\Y$ such that for each $\omega\in\Omega$ and $B\in \B(\Y)$, $\QM'(\omega,B)=\QM(\omega,\{x:g(\omega,x)\in B\})$.

The \emph{characteristic function} of a probability kernel $\QM$ on $\R^d$ is $\widehat\QM(\omega,t)\eqdef\int e^{\mathbf{i}t^{\top}x}\,\QM(\omega,dx)$, where $\mathbf{i}$ is the imaginary unit. Note that $\widehat\QM$ is a Carath\'eodory function on $\Omega\times \R^d$. For a random vector $X$ in $\R^d$ and a sub-$\sigma$-algebra $\F\subset\H$, the conditional characteristic function $\varphi_X(t\mid\F)\eqdef\E[\exp(\mathbf{i}t^{\top}X)\mid\F]$ satisfies $\varphi_X(t\mid\F)=\widehat\QM_X^{\F}(\csdot,t)$ a.s.\ for each $t\in\R^d$.

For an arbitrary family of random variables $\{X_{\alpha}:\alpha\in A\}$, $\esssup_{\alpha\in A}X_{\alpha}$ denotes its essential supremum, i.e., a random variable $X$ such that (a) $X\ge X_{\alpha}$ a.s.\ for each $\alpha\in A$; and (b) $X\le Z$ a.s.\ for any $Z$ with $Z\ge X_{\alpha}$ a.s.\ for each $\alpha\in A$. In particular, there is a countable $\seq{\alpha_k}\subset A$ with $X=\sup_k X_{\alpha_k}$ a.s.\ \citep[see, e.g.,][Section 6.5]{Chow/Teicher:97:Probability}. We refer to this as the \emph{countable-subfamily property}. The essential infimum is defined as $\essinf_{\alpha}X_{\alpha}\eqdef-\esssup_{\alpha\in A}(-X_{\alpha})$.

\section{Conditional Convergence}
\label{sec:cond_conv}

This section presents several modes of conditional convergence as well as the notion of asymptotic negligibility, which expresses the asymptotic tightness and equicontinuity of a sequence of probability kernels. These modes are useful for comparing two sequences of conditional laws. The comparison applies the bounded-Lipschitz metric of \cite{Dudley:68} to the merging framework of \cite{DAristotile/Diaconis/Freedman:88} and the conditional setting of \cite{Belyaev:00}. The following results are formulated in terms of probability kernels. The conditioning $\sigma$-algebras enter only through the formation of the kernels. Throughout this section, $\FS\eqdef\seq{\F_n}$ is an arbitrary sequence of sub-$\sigma$-algebras of $\H$.

\subsection{Modes of convergence}
\label{sec:modes}

We begin with conditional convergence in probability along a sequence of conditioning fields.

\begin{definition}
\label{def:CCP}
Let $\seq{X_n}$ be a sequence of random variables. We write $X_n\PRC{\FS}0$ a.s.\ if for any $\epsilon>0$, $\PR{\absin{X_n}>\epsilon\mid\F_n}\to 0$ a.s.\footnote{
	This implies convergence in probability by the dominated convergence theorem. The constant $\epsilon$ can also be replaced by a sequence of $\F_n$\-/measurable random variables $\{\nu_n\}$ with $\liminf_n\nu_n>0$ a.s.
}
\end{definition}

The primitive notion of closeness below is the Kolmogorov distance. Bounding this distance is the classical route to bootstrap consistency \citep[see, e.g.,][Chapter 3]{Shao:95:Bootstrap}. In the conditional form, it yields the validity of conditional quantiles in Section \ref{sec:cond_bootstrap}. For probability kernels $\QM$ and $\QM'$ on $\R$ it is given by
\[
	\DK{\QM,\QM'}\eqdef\sup_{x\in\R}\abs{\QM(\csdot,(-\infty,x])-\QM'(\csdot,(-\infty,x])},
\]
and, for random variables $X$ and $Z$ and sub-$\sigma$-algebras $\F\subset\G\subset\H$, we define
\[
	\DK{X,Z\mid\G,\F}\eqdef\DK{\QM_X^{\G},\QM_{Z}^{\F}},
\]
writing $\dK{X,Z\mid\F}$ when $\F=\G$.\footnote{
	The process $D_x\eqdef\abs{\QM(\csdot,(-\infty,x])-\QM'(\csdot,(-\infty,x])}$ is c\`adl\`ag in $x$, so the supremum can be taken over the rationals. Therefore, $\omega\mapsto \dK{\QM(\omega,\csdot),\QM'(\omega,\csdot)}$ is measurable with respect to any $\sigma$-algebra for which both kernels are measurable.
}
This distance is invariant under an $\F$\-/measurable increasing transformation applied to both arguments.

\begin{lemma}
\label{lemma:scale_inv}
Let $\F\subset\G\subset\H$ and let $g:\Omega\times\R\to\R$ be such that $g(\csdot,x)$ is $\F$\-/measurable for every $x\in\R$ and $g(\omega,\csdot)$ is strictly increasing and continuous for every $\omega$ outside a $\PM$-null set $N\in\F$. Then, for any random variables $X$ and $Z$,
\[
	\DK{g(\csdot,X),g(\csdot,Z)\mid\G,\F}=\DK{X,Z\mid\G,\F} \qtext{a.s.}
\]
In particular, if $a>0$ a.s., and $a$ and $b$ are $\F$\-/measurable, then $\dK{aX+b,aZ+b\mid\G,\F}=\dK{X,Z\mid\G,\F}$ a.s.
\end{lemma}

\begin{proof}
For $x\in\R$, define $\tau_x\eqdef\inf\{t\in\R:g(\csdot,t)\ge x\}\in\R\cup\{\pm\infty\}$ on $N^{c}$, and $\tau_x\eqdef 0$ on $N$. Then we have,
\[
	\{\tau_x\le t\}=(\{g(\csdot,t)\ge x\}\cap N^{c})\cup(N\cap\{t\ge0\})\in\F,
\]
so that $\tau_x$ is $\F$\-/measurable, and $\{g(\csdot,X)\le x\}=\{X\le\tau_x\}$ on $N^c$. Evaluating the conditional cdf at $\tau_x$ yields $\PR{X\le\tau_x\mid\G}=\FD{X}{\G}(\csdot,\tau_x)$ a.s., where we set $\FD{X}{\G}(\csdot,-\infty)\eqdef 0$ and $\FD{X}{\G}(\csdot,\infty)\eqdef 1$. Hence, $\FD{g(\csdot,X)}{\G}(\csdot,x)=\FD{X}{\G}(\csdot,\tau_x)$ for all $x\in\R$ a.s.\ by the right-continuity in $x$ of both sides; likewise for $Z$ with $\F$ in place of $\G$. Outside $N$, the map $g(\omega,\csdot)$ is a strictly increasing continuous bijection onto its range $(c,d)$, and $x\mapsto\tau_x$ is its inverse on $(c,d)$. Thus, $\tau_x$ attains every real value as $x$ ranges over $(c,d)$. For $x\notin(c,d)$ one has $\tau_x=\pm\infty$ and the two conditional cdfs agree. Therefore, the suprema coincide.
\end{proof}

For a probability kernel $\QM$ on $\R$, the \emph{L\'evy concentration function} is given by
\[
	\LC{\epsilon,\QM}\eqdef\sup_{x\in\R}\QM(\csdot,[x,x+\epsilon]).
\]
The largest atom of $\QM$ is recovered as
\[
	\LCZ{\QM}\eqdef\lim_{\epsilon\searrow 0}\LC{\epsilon,\QM}=\sup_{x\in\R}\QM(\csdot,\{x\}).
\]
In particular, $\lCZ{\QM}=0$ when the cdf of $\QM$ is continuous.\footnote{
	The L\'evy concentration function can be written as
	\begin{align*}
		\LC{\epsilon,\QM}&=\inf_{\eta\in \Q\cap(0,\infty)}\sup_{q\in \Q}\QM(\csdot,[q,q+\epsilon+\eta]),
	\end{align*}
	and thus $\omega\mapsto\LC{\epsilon,\QM(\omega)}$ is measurable. Moreover, $\omega\mapsto\lCZ{\QM(\omega)}$ is measurable as a decreasing limit of measurable functions, and $\lCZ{\QM}\le\lC{\epsilon,\QM}$ for every $\epsilon>0$.
}
For a random variable $X$ and $\F\subset\H$, we abbreviate $\lC{\epsilon,X\mid\F}\eqdef\lC{\epsilon,\QM_X^{\F}}$ and $\lCZ{X\mid \F}\eqdef \lCZ{\QM_X^{\F}}$.

We introduce the notion of an asymptotically negligible random function, which is used to state the conditional versions of asymptotic tightness and equicontinuity.

\begin{definition}
\label{def:asy_negligible}
Let $A$ be a set and let $\seq{f_n}$ be a sequence of functions $f_n:\Omega\times A\to [0,\infty]$ such that $f_n(\csdot,\alpha)$ is $\H$\-/measurable for each $\alpha\in A$. The sequence $\seq{f_n}$ is \textit{asymptotically negligible} (a.n.) if
\[
	\essinf_{\alpha\in A}\limsup_{n\to\infty}f_n(\csdot,\alpha)=0 \qtext{a.s.}
\]
In particular:
\begin{enumerate}
	\item[--] A sequence of probability kernels $\seq{\QM_n}$ on $\X$ is \textit{asymptotically tight} if the family $\seq{\QM_n(\csdot,K^{c}),K\in\K_{\X}}$ is a.n.

	\item[--] A sequence of probability kernels $\seq{\QM_n}$ on $\R$ is \emph{asymptotically equicontinuous} if the family $\seq{\lC{\csdot,\QM_n}}$ is a.n.
\end{enumerate}

A sequence of $\X$\-/valued random elements $\seq{X_n}$ is $\FS$\-/\textit{asymptotically tight} if $\seq{\QM_{X_n}^{\F_n}}$ is asymptotically tight, and a sequence of random variables $\seq{X_n}$ is \emph{$\FS$\-/asymptotically equicontinuous} if $\seq{\QM_{X_n}^{\F_n}}$ is asymptotically equicontinuous.\footnote{
	Since $\sup_{x}(F(x+\epsilon)-F(x))\le\lC{\epsilon,\QM}\le\sup_{x}(F(x+2\epsilon)-F(x))$ for the cdf $F$ of $\QM$, this is the asymptotic uniform equicontinuity of the cdfs, holding $\omega$-wise outside a null set.
}
\end{definition}

By the countable-subfamily property, there exists a countable $\seq{\alpha_k}\subset A$ with
\begin{equation}
	\label{eq:countable_equiv}
	\essinf_{\alpha\in A}\limsup_{n\to\infty} f_n(\csdot,\alpha)=\inf_k\limsup_{n\to\infty}f_n(\csdot,\alpha_k) \qtext{a.s.},
\end{equation}
so that asymptotic negligibility implies that, for a fixed countable family of indices, the right-hand side of \eqref{eq:countable_equiv} is zero, outside a single null set. For example, the following result shows that the asymptotic tightness of a sequence of kernels is equivalent to uniform tightness in the usual sense, holding $\omega$-wise outside a null set.

\begin{lemma}
\label{lemma:asy_tight_equiv}
Let $\seq{\QM_n}$ be a sequence of probability kernels on $\X$. Then, the following statements are equivalent.
\begin{enumerate}[label=\tn{(\alph*)}]
	\item $\seq{\QM_n}$ is asymptotically tight.

	\item For almost every $\omega\in\Omega$, the sequence $\seq{\QM_n(\omega,\csdot)}$ is uniformly tight, i.e., for every $\epsilon>0$, there exists $K\in\K_{\X}$, depending on $\omega$ and $\epsilon$, with $\sup_{n\ge1}\QM_n(\omega,K^{c})\le\epsilon$.

	\item There exists a nondecreasing sequence $\{K_m\}$ in $\K_{\X}$ such that $\sup_{n\ge 1}\QM_n(\csdot,K_m^{c})\to0$ a.s.\ as $m\to\infty$.
\end{enumerate}
\end{lemma}

\begin{proof}
The implication (a) $\Rightarrow$ (b) follows from the countable-subfamily property and Ulam's theorem \citep[Theorem 7.1.4]{Dudley:02:RAP}, and (c) $\Rightarrow$ (a) holds because for each $m\ge 1$,
\[
	\essinf_{K\in\K_{\X}}\limsup_{n\to\infty}\QM_n(\csdot,K^{c})\le\sup_{n\ge 1}\QM_n(\csdot,K_m^{c}) \qtext{a.s.},
\]
and the right-hand side converges to $0$ a.s.\ as $m\to\infty$.

\smallskip

\noindent (b) $\Rightarrow$ (c) Let $\seq{x_1,x_2,\ldots}$ be a countable dense subset of $\X$, and for each $j,k\in\N$, let $C_{j,k}\eqdef\bigcup_{i\le k} \overline{B(x_i,1/j)}$, where $B(x,r)\eqdef \{z\in\X:\varrho_{\X}(z,x)<r\}$. Fix $\epsilon > 0$ and define
\[
	\kappa_{j}(\omega)\eqdef\inf\left\{k\ge1:\sup_{n\ge 1}\QM_n(\omega, C_{j,k}^{c})\le\epsilon2^{-j}\right\}\in\N\cup\{\infty\}.
\]
Note that the random variable $\kappa_{j}$ is a.s.\ finite. Indeed, fix $\omega$ for which $\seq{\QM_n(\omega,\csdot)}$ is uniformly tight. Then there exists $K\in\K_{\X}$ such that $\sup_{n\ge 1}\QM_n(\omega,K^c)\le \epsilon 2^{-j}$. Since $K$ is covered by $\bigcup_{i\ge 1}B(x_i,1/j)$, there exists a finite subcover $\bigcup_{i\le k^*}B(x_i,1/j)\subset C_{j, k^*}$ of $K$ and $\kappa_{j}(\omega)\le k^*<\infty$.

For each $j\ge 1$, let $k_{j}>0$ be the smallest $k>0$ such that $\PR{\kappa_{j}>k}\le\epsilon 2^{-j}$ and set $K\eqdef\bigcap_{j\ge 1}C_{j,k_{j}}$. Clearly, $K$ is compact, and on the event $G\eqdef\bigcap_{j\ge 1}\{\kappa_{j}\le k_{j}\}$, having $\PM(G^c)\le \epsilon$, for each $n\ge 1$,
\[
	\QM_n(\csdot,K^{c})\le\sum_{j\ge 1}\QM_n(\csdot,C_{j,k_j}^{c})\le\epsilon,
\]
so that $\sup_n\QM_n(\csdot,K^{c})\le\epsilon$ on $G$.

Consequently, for each $m\ge 1$, there is $K'_m\in\K_{\X}$ with $\PR{\sup_n\QM_n(\csdot,K'^c_m)>1/m}\le 1/m$. Let $K_m\eqdef \bigcup_{i\le m}K'_i$. Since $\seq{K_m}$ is nondecreasing, $\xi_m\eqdef\sup_n\QM_n(\csdot,K_m^{c})$ is nonincreasing in $m$ and therefore converges to a limit $\xi_{\infty}$. Since for any $\delta> 0$ and $m\ge 1/\delta$, $\PR{\xi_{\infty}>\delta}\le\PR{\xi_m>1/m}\le 1/m$, we have $\xi_{\infty}=0$ a.s.
\end{proof}

\begin{remark}
\label{rem:essinf}
(1) The essential infimum depends only on the equivalence classes of $f_n(\csdot,\alpha)$, so that asymptotic negligibility in Definition \ref{def:asy_negligible} is unaffected by the choice of the versions representing the kernels. The countable family in \eqref{eq:countable_equiv} is not explicitly exhibited in general, but one can take it to be $\Q \cap (0,\infty)$ when $A=(0,\infty)$ and the map $\alpha\mapsto f_n(\omega,\alpha)$ is monotone, e.g., the L\'evy concentration function.

(2) Tightness is preserved by continuous maps. Specifically, suppose that $\seq{X_n}$ is $\FS$\-/asymptotically tight, each $X_n\in C$ a.s.\ for some closed $C\subset \X$, and $f:C\to \Y$ is continuous. Then, for $K\in\K_{\X}$, both $C\cap K$ and $f(C\cap K)$ are compact, and
\begin{align*}
	&\essinf_{K\in \K_{\X}}\limsup_{n\to\infty}\PR{f(X_n)\notin f(C\cap K)\mid\F_n} \\
	&\qquad\le\essinf_{K\in \K_{\X}}\limsup_{n\to\infty}\PR{X_n\notin K\mid\F_n}=0 \qtext{a.s.}
\end{align*}
by the countable-subfamily property.
That is, $\seq{f(X_n)}$ is $\FS$\-/asymptotically tight.
\end{remark}

\subsection{Merging and \texorpdfstring{$\FS$}{F}-weak convergence}
\label{sec:merging}

Given a Polish space $\X$ and probability measures $\mu$ and $\nu$ on $\X$, let $\BL{\X}\eqdef \{f:\X\to\R: \norm{f}_{\infty}\le 1, \Lip(f)\le 1\}$, and let
\[
	\DBL{\mu,\nu}\eqdef\sup\left\{\abs{\int f\,d\mu-\int f\,d\nu}:f\in \BL{\X}\right\}
\]
denote the bounded-Lipschitz metric on the space of probability measures on $(\X,\B(\X))$.\footnote{
	For $\F$\-/measurable probability kernels $\QM$ and $\QM'$, the map $\omega\mapsto \dBL{\QM(\omega,\csdot),\QM'(\omega,\csdot)}$ is $\F$\-/measurable. To see this, let $\P$ be the set of probability measures on $\X$ with the weak topology and its Borel $\sigma$-algebra. Since $\BLM$ metrizes this topology \citep[\S11.3]{Dudley:02:RAP}, it is continuous on $\P\times\P$, and thus Borel. The assertion follows because $\P$ is separable and $\omega \mapsto (\QM(\omega,\csdot), \QM'(\omega,\csdot))$ is $(\F,\B(\P)\otimes\B(\P))$\-/measurable.
}

\begin{definition}
\label{def:CWC}
Let $\seq{\X_n}$ be a sequence of Polish spaces and for each $n\ge 1$, let $\QM_n$ and $\QM_n'$ be probability kernels on $\X_n$. We say that $\seq{\QM_n}$ \emph{merges} with $\seq{\QM_n'}$, written as $\QM_n\mrg \QM_n'$, if $\dBL{\QM_n,\QM_n'}\to 0$ a.s. Sequences of random elements $\seq{X_n}$ and $\seq{Z_n}$, where for each $n\ge 1$, $X_n$ and $Z_n$ are $\X_n$-valued, \emph{merge $\FS$\-/weakly} if $\QM_{X_n}^{\F_n}\mrg \QM_{Z_n}^{\F_n}$.
\end{definition}

When $\X_n\equiv\X$, since $\BLM$ metrizes weak convergence, merging with a constant kernel $\QM_n'\equiv\QM$ reduces to the weak convergence of $\QM_n(\omega,\csdot)$ to $\QM(\omega,\csdot)$ for almost every $\omega\in\Omega$, written as $\QM_n\Rightarrow\QM$. For random elements we call this \emph{$\FS$-weak convergence} and write $X_n\to X$, $\FS$\-/weakly w.r.t.\ $\F$, when $\QM_n=\QM_{X_n}^{\F_n}$ and $\QM$ is the conditional law of $X$ given $\F\subset\F_{\infty}$ (see Lemma \ref{lemma:tail_join} below). \emph{$\F$-weak convergence}, written $X_n\to X$, $\F$\-/weakly, is the special case of a constant sequence $\F_n\equiv\F$, in which $\QM=\QM_X^{\F}$.

\begin{remark}
(1) The terminology follows the literature on \emph{merging} of probability measures; see \cite{DAristotile/Diaconis/Freedman:88} for the notions compatible with the weak topology and the relations among them and \cite{Davydov/Rotar:09} for a systematic treatment. The version for random distributions closest to Definition \ref{def:CWC} is the almost sure form of the weakly approaching sequences of \citet[Definitions 2 and 3]{Belyaev:00}, who require $\int h\,d(\QM_n-\QM_n')\to 0$ for each $h\in\Cb{\X}$, rather than convergence in the bounded-Lipschitz metric. For nonrandom measures, these definitions agree under tightness \citep[Remark 2]{Davydov/Rotar:09} and may differ otherwise \citep[Example 1.1 and Theorem 2.8]{DAristotile/Diaconis/Freedman:88}.

(2) $\FS$\-/weak convergence of $\{X_n\}$ to $\QM$ implies its ordinary weak convergence to the mixture $\E[\QM]$ because, for $f\in\Cb{\X}$,
\[
	\E[f(X_n)]=\E[\E[f(X_n)\mid\F_n]]\to \int \left(\int f\,d\QM \right) d\PM
\]
by the dominated convergence theorem. Likewise, $\FS$\-/weak convergence implies the convergence of $\QM_{X_n}^{\F_n}$ to $\QM$ in the sense of Markov kernels \citep[see][Section 1.13]{vanderVaart/Wellner:23}. For a nonrandom limit $L$, $\FS$\-/weak convergence is the almost sure version of the conditional weak convergence of \citet[Section 1.13]{vanderVaart/Wellner:23}. Finally, for a fixed $\F\subset\H$, $\F$\-/weak convergence $X_n\to X$ implies $\F$\-/stable convergence to $\QM_X^{\F}$ because $\E[Zf(X_n)]=\E[Z\int f\,d\QM_{X_n}^{\F}]\to\E[Z\int f\,d\QM_X^{\F}]$ for every $f\in\Cb{\X}$ and $Z\in L^1(\F)$ \citep[see, e.g.,][]{Hausler/Luschgy:15}, and it is strictly stronger.
\end{remark}

The properties of merging summarized below mirror the classical weak-convergence toolbox with one exception. When the comparison sequence changes with $n$, there is no fixed object providing tightness, as the continuity of the limit characteristic function at the origin does in the classical case, and therefore asymptotic tightness must be assumed.

For a function $g:\X\to\Y$ and $\delta>0$, let
\[
	\varpi_g(\delta)\eqdef\sup\left\{\varrho_{\Y}(g(x),g(z)):x,z\in\X,\ \varrho_{\X}(x,z)\le\delta\right\}
\]
denote its modulus of continuity. For a Carath\'eodory function $f$, we write $\varpi_f(\delta)$ for $\omega\mapsto\varpi_{f(\omega,\csdot)}(\delta)$.\footnote{
	Since $(x,z)\mapsto\varrho_{\Y}(f(\omega,x),f(\omega,z))$ is continuous for each $\omega\in\Omega$, the supremum may be taken over a countable dense subset of $\{(x,z)\in\X\times\X:\varrho_{\X}(x,z)\le\delta\}$, so that $\varpi_f(\delta)$ is $\H$\-/measurable.
}

\begin{lemma}
\label{lemma:merging}
Let $\seq{\X_n}$ be a sequence of Polish spaces and for each $n\ge 1$, let $\QM_n$ and $\QM_n'$ be probability kernels on $\X_n$.
\begin{enumerate}[label=\tn{(\alph*)}]
	\item \tn{(Merging together)} Suppose that for each $n\ge 1$, $X_{n}$ and $Z_{n}$ are $\X_n$-valued random elements. If $\QM_{X_n}^{\F_n}\mrg \QM_n$ and $\varrho_{\X_n}(X_n,Z_n)\PRC{\FS}0$ a.s., then $\QM_{Z_n}^{\F_n}\mrg \QM_n$.

	\item \tn{(Continuous mapping)} Suppose that $\QM_n\mrg \QM_n'$. Then $\QM_n\circ g_n^{-1}\mrg\QM_n'\circ g_n^{-1}$ if either \tn{(i)} for each $n\ge 1$, $g_n:\Omega\times\X_n\to\Y_n$ is a Carath\'eodory function into a Polish space $\Y_n$, and $\seq{\varpi_{g_n}(\delta),\delta>0}$ is a.n., or \tn{(ii)} $\X_n\equiv\X$, $\Y_n\equiv\Y$, $g_n\equiv g$ for a continuous function $g:\X\to\Y$, and $\seq{\QM_n}$ is asymptotically tight.

	\item \tn{(Cram\'er--Wold)} Let $\X_n\equiv\R^d$ and let $D$ be a countable and dense subset of $\R^d$. Suppose that $\seq{\QM_n}$ is asymptotically tight. Then $\QM_n\mrg \QM_n'$ if and only if $\QM_n\circ(t^{\top}\csdot)^{-1}\mrg \QM_n'\circ(t^{\top}\csdot)^{-1}$ for each $t\in D$.

	\item \tn{(L\'evy continuity 1)} Let $\X_n\equiv\R^d$ and let $D$ be a countable and dense subset of $\R^d$. Suppose that $\seq{\QM_n}$ and $\seq{\QM_n'}$ are asymptotically tight. Then $\QM_n\mrg \QM_n'$ if and only if $\widehat\QM_n(\csdot,t)-\widehat\QM_n'(\csdot,t)\to0$ a.s.\ for each $t\in D$.

	\item \tn{(L\'evy continuity 2)} Let $\X_n\equiv\R^d$ and suppose that $\seq{\QM_n}$ is asymptotically tight. Then $\QM_n\mrg\QM_n'$ if and only if $\widehat\QM_n(\csdot,t)-\widehat\QM_n'(\csdot,t)\to0$ a.s.\ for each $t\in\R^d$.
\end{enumerate}
\end{lemma}

\begin{remark}
\label{rem:merging1}
Since $\QM_{g_n(\csdot,X_n)}^{\F_n}=\QM_{X_n}^{\F_n}\circ g_n^{-1}$ a.s.\ for an $\F_n\otimes\B(\X_n)$\-/measurable function $g_n:\Omega\times\X_n\to\Y_n$ \citep[see, e.g.,][Theorem 8.5]{Kallenberg:21:FMP}, parts (b) and (c) of Lemma \ref{lemma:merging} can be formulated in terms of sequences of random elements, i.e., using $\seq{g_n(\csdot,X_n)}$ in (b)(i), $\seq{g(X_n)}$ in (b)(ii), and $\seq{t^{\top}X_n}$ in (c).
\end{remark}

\begin{proof}[Proof of Lemma \ref{lemma:merging}]
(a) Fix $\epsilon>0$. For $f_n\in \BL{\X_n}$, we have $\abs{f_n(x)-f_n(z)}\le 2\wedge \varrho_{\X_n}(x,z)$, so that for almost every $\omega\in\Omega$,
\begin{align*}
	\DBL{\QM_{X_n}^{\F_n}(\omega,\csdot),\QM_{Z_n}^{\F_n}(\omega,\csdot)}&\le\int (2\wedge \varrho_{\X_n}(x,z))\,\QM_{(X_n,Z_n)}^{\F_n}(\omega,d(x,z)) \\
	&\le\epsilon+2\QM_{(X_n,Z_n)}^{\F_n}(\omega,\{\varrho_{\X_n}(x,z)>\epsilon\}).
\end{align*}
Since $\QM_{(X_n,Z_n)}^{\F_n}(\csdot,\{\varrho_{\X_n}(x,z)>\epsilon\})\to 0$ a.s., by the triangle inequality,
\begin{align*}
	\limsup_{n\to\infty}\DBL{\QM_{Z_n}^{\F_n},\QM_n}&\le\limsup_{n\to\infty}\DBL{\QM_{Z_n}^{\F_n},\QM_{X_n}^{\F_n}} \\
	&\quad+\limsup_{n\to\infty}\DBL{\QM_{X_n}^{\F_n},\QM_n}\le\epsilon \qtext{a.s.}
\end{align*}
Taking $\epsilon$ along a sequence converging to zero implies the result.

\smallskip

\noindent(b) Assume (i) and fix $\omega$ outside the null sets of the merging hypothesis and of the asymptotic negligibility of the moduli. Take $f_n\in \BL{\Y_n}$ and define $h_n\eqdef f_n\circ g_n(\omega,\csdot)$, so that $\norm{h_n}_{\infty}\le 1$ and $\varpi_{h_n}(\delta)\le 2\wedge\varpi_{g_n(\omega,\csdot)}(\delta)$ for each $\delta>0$. For $m\ge 1$, the Pasch--Hausdorff envelope
\[
	h_{n,m}(x)\eqdef\inf_{z\in\X_n}\left[h_n(z)+m\varrho_{\X_n}(x,z)\right]
\]
is $m$\-/Lipschitz with $\norm{h_{n,m}}_{\infty}\le\norm{h_n}_{\infty}$ \citep[see, e.g.,][Theorem 2.1]{Gutev:20} and
\begin{align*}
	\norm{h_n-h_{n,m}}_{\infty}&=\sup_{t\ge 0}\left[\varpi_{h_n}(t)-mt\right] \\
	&\le\sup_{t\ge 0}\left[(2\wedge\varpi_{g_n(\omega,\csdot)}(t))-mt\right]\eqqcolon \epsilon_{n,m}.
\end{align*}
Therefore, using the triangle inequality, for each $m\ge 1$,
\begin{align}
\label{eq:bl_image_bound}
	\begin{aligned}
		&\DBL{\QM_n(\omega,\csdot)\circ g_n(\omega,\csdot)^{-1},\QM_n'(\omega,\csdot)\circ g_n(\omega,\csdot)^{-1}} \\
		&\qquad\le 2\epsilon_{n,m}+m \cdot\DBL{\QM_n(\omega,\csdot),\QM_n'(\omega,\csdot)}.
	\end{aligned}
\end{align}
The expression inside the supremum in the definition of $\epsilon_{n,m}$ is at most $2\wedge\varpi_{g_n(\omega,\csdot)}(2/m)$ for $t\le2/m$ and negative for $t>2/m$. Thus, $\limsup_n\epsilon_{n,m}\le \limsup_{n}(2\wedge\varpi_{g_n(\omega,\csdot)}(2/m))\to 0$ as $m\to\infty$. The result follows by sending $n\to\infty$ and then $m\to\infty$ in \eqref{eq:bl_image_bound}.

\smallskip

For alternative (ii) of (b) and for parts (c)--(d), we argue $\omega$-wise along subsequences. Fix $\omega$ outside the countably many null sets arising from the stated conditions. By Lemma \ref{lemma:asy_tight_equiv}, the sequence $\seq{\QM_n(\omega,\csdot)}$ is uniformly tight, so that every subsequence has a further subsequence converging weakly by Prokhorov's theorem.

\smallskip

\noindent(b) Assume (ii) and suppose the conclusion fails at $\omega$, so that there are $\epsilon>0$ and a subsequence, both possibly depending on $\omega$, along which $\dBL{\QM_n(\omega,\csdot)\circ g^{-1},\QM_n'(\omega,\csdot)\circ g^{-1}}\ge\epsilon$. Along a further subsequence $\QM_{n_k}(\omega,\csdot)\Rightarrow\mu$, and since
\[
	\DBL{\QM_{n_k}'(\omega,\csdot),\mu}\le\DBL{\QM_{n_k}'(\omega,\csdot),\QM_{n_k}(\omega,\csdot)}+\DBL{\QM_{n_k}(\omega,\csdot),\mu}\to0,
\]
we have $\QM_{n_k}'(\omega,\csdot)\Rightarrow\mu$. The continuous mapping theorem yields $\QM_{n_k}(\omega,\csdot)\circ g^{-1}\Rightarrow\mu\circ g^{-1}$ and $\QM_{n_k}'(\omega,\csdot)\circ g^{-1}\Rightarrow\mu\circ g^{-1}$, so their $\BLM$\-/distance tends to $0$, a contradiction.

\smallskip

\noindent(c) Necessity follows from (b)(i) with the uniformly continuous $g(x)=t^{\top}x$. For sufficiency, write $\pi_t(x)\eqdef t^{\top}x$. First, we show the transfer of tightness for merging sequences of probability measures on $\R$.

\begin{claim*}
	Let $\seq{\mu_n}$ and $\seq{\nu_n}$ be sequences of probability measures on $\R$ such that $\seq{\mu_n}$ is uniformly tight and $\mu_n\mrg\nu_n$. Then, $\seq{\nu_n}$ is uniformly tight.
\end{claim*}

\begin{subproof}
For $y>0$, let $f_y$ be the piecewise linear function equal to $1$ on $\{\absin{x}\le y\}$, to $0$ on $\{\absin{x}> y+1\}$, and linear in $\absin{x}$ in between, so that $f_y\in\BL{\R}$. Since $\ind\{\absin{x}\le y\}\le f_y\le\ind\{\absin{x}\le y+1\}$, we have
\[
	\nu_n(\{\abs{x}>y+1\})\le \mu_n(\{\abs{x}>y\})+\DBL{\mu_n,\nu_n}.
\]
As $\DBL{\mu_n,\nu_n}\to 0$ as $n\to\infty$, $\seq{\nu_n}$ is asymptotically tight, and thus uniformly tight.
\end{subproof}

Since $\absin{t^{\top}x}\le\norm{t}\norm{x}$, the tightness of $\seq{\QM_n(\omega,\csdot)}$ transfers to $\seq{\QM_n(\omega,\csdot)\circ\pi_t^{-1}}$, and the preceding claim carries it to $\seq{\QM_n'(\omega,\csdot)\circ\pi_t^{-1}}$ for every $t\in D$. Let $T$ denote the matrix whose rows are transposed linearly independent vectors $t_1,\ldots,t_d\in D$ and $\kappa\eqdef\sqrt{d}\,\norm{T^{-1}}_{\mathrm{op}}$. Then $\norm{x}\le\kappa\max_{1\le j\le d}\absin{t_j^{\top}x}$, so that $\{\norm{x}>y\}\subset\bigcup_{1\le j\le d}\{\absin{t_j^{\top}x}>y/\kappa\}$ and
\[
	\QM_n'(\omega,\{\norm{x}>y\})\le\sum_{1\le j\le d}(\QM_n'(\omega,\csdot)\circ\pi_{t_j}^{-1})(\{\abs{u}>y/\kappa\}).
\]
Therefore, $\seq{\QM_n'(\omega,\csdot)}$ is uniformly tight. Suppose the conclusion fails at $\omega$, so that along a subsequence, $\dBL{\QM_n(\omega,\csdot),\QM_n'(\omega,\csdot)}\ge\epsilon$ for some $\epsilon>0$. Extract a further subsequence such that $\QM_{n_k}(\omega,\csdot)\Rightarrow\mu$ and $\QM_{n_k}'(\omega,\csdot)\Rightarrow\nu$, where $\mu\ne\nu$. For $t\in D$, the continuous mapping theorem implies $\mu\circ\pi_t^{-1}=\nu\circ\pi_t^{-1}$, whence $\widehat\mu(t)=\widehat\nu(t)$. Since characteristic functions are continuous and $D$ is dense, $\mu=\nu$. In this case, however, $\dBL{\QM_{n_k}(\omega,\csdot),\QM_{n_k}'(\omega,\csdot)}\to 0$, a contradiction.

\smallskip

\noindent(d) We show necessity first. For $t\ne0$, the functions $x\mapsto\cos(t^{\top}x)$ and $x\mapsto\sin(t^{\top}x)$ are bounded by $1$ and Lipschitz with constant $\norm{t}$. Therefore,
\[
	\absin{\widehat\QM_n(\omega,t)-\widehat\QM_n'(\omega,t)}\le2(1\vee\norm{t})\DBL{\QM_n(\omega,\csdot),\QM_n'(\omega,\csdot)}\to 0
\]
as $n\to\infty$. As for sufficiency, if the conclusion fails along a subsequence, extract a further subsequence such that $\QM_{n_k}(\omega,\csdot)\Rightarrow\mu$ and $\QM_{n_k}'(\omega,\csdot)\Rightarrow\nu$, where $\mu\ne\nu$. Then $\widehat\mu(t)=\widehat\nu(t)$ for all $t\in D$, and both characteristic functions are continuous, so that they agree everywhere and $\mu=\nu$ by uniqueness, a contradiction.

\smallskip

\noindent(e) Necessity follows from the bound in the proof of (d). For sufficiency, let $\lambda$ denote the Lebesgue measure on $\R^d$. Since the map $(\omega,t)\mapsto\widehat\QM_n(\omega,t)-\widehat\QM_n'(\omega,t)$ is jointly measurable, by Fubini's theorem there is a null set outside of which it converges to zero for $\lambda$\-/almost every $t$. Fix $\omega$ outside this set at which $\seq{\QM_n(\omega,\csdot)}$ is uniformly tight.

\begin{claim*}
	Let $\seq{\mu_n}$ and $\seq{\nu_n}$ be sequences of probability measures on $\R^d$ such that $\seq{\mu_n}$ is uniformly tight and $\widehat\mu_n(t)-\widehat\nu_n(t)\to0$ for $\lambda$\-/almost every $t$. Then $\mu_n\mrg\nu_n$.
\end{claim*}

\begin{subproof}
For $\sigma>0$, let $\gamma_{\sigma}$ denote a Gaussian law on $\R^d$ with mean zero and covariance matrix $\sigma^2 I_d$. The convolved measures $\mu_n*\gamma_{\sigma}$ and $\nu_n*\gamma_{\sigma}$ have continuous densities $p_n$ and $p_n'$. Thus, by the inversion formula and the dominated convergence theorem, we get
\[
	\eta_n\eqdef\sup_{x\in\R^d}\abs{p_n(x)-p_n'(x)}\le \frac{1}{(2\pi)^d}\int\absin{\widehat\mu_n(t)-\widehat\nu_n(t)}e^{-\frac{\sigma^2\norm{t}^2}{2}}\,\lambda(dt)\to 0
\]
as $n\to\infty$. Fix $\epsilon>0$ and $K\in \K_{\R^d}$ with $(\mu_n*\gamma_{\sigma})(K^c)\le\epsilon$ for every $n\ge 1$. Since
\[
	\int_{K^c}p_n'\,d\lambda \le\int_{K^c}p_n\,d\lambda+\lambda(K)\eta_n,
\]
we have
\[
	\int\abs{p_n-p_n'}\,d\lambda\le \lambda(K)\eta_n + \int_{K^c} (p_n + p_n')\, d\lambda \le 2\lambda(K)\eta_n+2\epsilon.
\]
Note that $\dBL{\mu_n*\gamma_{\sigma},\nu_n*\gamma_{\sigma}}\le\int\abs{p_n-p_n'}\,d\lambda$, and for a probability measure $\mu$ on $\R^d$, $\dBL{\mu,\mu*\gamma_{\sigma}}\le\sigma\sqrt{d}$. Thus, $\limsup_n\dBL{\mu_n,\nu_n}\le2\sigma\sqrt{d}+2\epsilon$, by the triangle inequality, and, consequently, the arbitrariness of $\sigma$ and $\epsilon$ implies the result.
\end{subproof}

The result follows from the claim above with $\mu_n=\QM_n(\omega,\csdot)$ and $\nu_n=\QM_n'(\omega,\csdot)$.
\end{proof}

If $g_n(\omega,\csdot)$ in part (b)(i) is Lipschitz with constant $L_n(\omega)$, then for $f\in\BL{\Y_n}$, the function $f\circ g_n(\omega,\csdot)/(1\vee L_n(\omega))\in\BL{\X_n}$. In this case,
\begin{align}
	\label{eq:lip_image}
	\begin{aligned}
		&\DBL{\QM_n(\omega,\csdot)\circ g_n(\omega,\csdot)^{-1},\QM_n'(\omega,\csdot)\circ g_n(\omega,\csdot)^{-1}} \\
		&\qquad\le(1\vee L_n(\omega))\,\DBL{\QM_n(\omega,\csdot),\QM_n'(\omega,\csdot)}.
	\end{aligned}
\end{align}
In particular, (b)(i) holds if $\limsup_nL_n<\infty$ a.s.

\begin{remark}
\label{rem:merging2}
(1) Part (a) of Lemma \ref{lemma:merging} is the conditional form of Theorem 2(A) of \cite{Davydov/Rotar:09}. Theorems 1 and 2 of \cite{Belyaev:00} are unconditional and in-probability counterparts of part (e). When $\X_n$ is fixed and $g$ is nonrandom, part (b)(i) is Corollary 6 of \cite{Davydov/Rotar:09}, whose hypothesis follows from Definition \ref{def:CWC} via their Theorem 3. The argument above is given to keep the passage elementary. The transfer of tightness for sequences of probability measures on $\R$ in the proof of part (c) is a special case of Theorem 3 and Remark 2 of \cite{Davydov/Rotar:09}; see also Lemmas 5 and 6 of \cite{Belyaev:00}.

(2) A constant comparison $\QM_n'\equiv\QM$ is tight by Ulam's theorem. Since merging is symmetric, part (e) implies that $\QM_n\Rightarrow\QM$ if and only if $\widehat\QM_n(\csdot,t)\to\widehat\QM(\csdot,t)$ a.s.\ for each $t\in\R^d$. In particular, $X_n\to X$, $\FS$\-/weakly w.r.t.\ $\F$, if and only if $\varphi_{X_n}(t\mid\F_n)\to\varphi_X(t\mid\F)$ a.s.\ for each $t$. Sufficiency is Theorem 2.6 of \cite{Berti:06}, and Theorem 1.13.1(iv) of \cite{vanderVaart/Wellner:23} is the counterpart for a nonrandom limit and convergence in probability.
\end{remark}

Example \ref{ex:cont_map} below shows that the conclusion of Lemma \ref{lemma:merging}(b) under alternative (ii) fails without tightness. By Theorem 4.4 of \cite{DAristotile/Diaconis/Freedman:88}, a map preserves merging in the bounded-Lipschitz metric if and only if it is uniformly continuous, whereas merging against all bounded continuous functions is preserved by every continuous map (see their Theorem 4.2). Alternative (ii) imposes only continuity and restricts $\seq{\QM_n}$ instead, since a tight sequence concentrates on compact sets, on which a continuous map is uniformly continuous. The function $g$ in the example also shows that, without tightness, merging against all bounded continuous functions is strictly stronger than merging in the bounded-Lipschitz metric.

\begin{example}
\label{ex:cont_map}
Let $\delta_x$ denote the Dirac measure at $x$. Take $\F_n\equiv\H_{\varnothing}$, $\QM_n\eqdef\delta_{2n}$, and $\QM_n'\eqdef\delta_{2n+1/n}$, so that $\dBL{\QM_n,\QM_n'}=1/n\to0$ while neither sequence is tight \citep[see][Example 1.1]{DAristotile/Diaconis/Freedman:88}. Define $g:\R\to[0,1]$ to be the tent function on each $I_k\eqdef[2k,2k+2/k]$, $k\ge1$, that vanishes at the two endpoints, equals $1$ at the midpoint $2k+1/k$, linear on each half, and $0$ outside $\bigcup_{k\ge1}I_k$. The intervals $I_k$ have disjoint interiors and the only point lying in two of them is $4$, an endpoint of both $I_1$ and $I_2$, where the two tents vanish, so $g$ is well defined and continuous. It is bounded but not uniformly continuous, and taking $k=n$ yields $g(2n)=0$ and $g(2n+1/n)=1$, implying $\QM_n\circ g^{-1}=\delta_0$ and $\QM_n'\circ g^{-1}=\delta_1$ for every $n$. Consequently, the images do not merge.
\end{example}

For a sequence of kernels that converges to a limit, the next result identifies the $\sigma$-algebra on which the limit must live.

\begin{lemma}
\label{lemma:tail_join}
Let $\seq{\QM_n}$ be a sequence of $\F_n$\-/measurable probability kernels on $\X$. If $\QM_n\Rightarrow\QM$ for a probability kernel $\QM$ on $\X$, then $\QM$ agrees almost surely with an $\F_{\infty}$\-/measurable kernel. In particular, the following hold.
\begin{enumerate}[label=\tn{(\alph*)}]
	\item If $\seq{\F_n}$ is nondecreasing, then $\F_{\infty}=\F^{\vee}$.
	\item If $\seq{\F_n}$ is nonincreasing, then $\F_{\infty}=\F^{\wedge}$.
	\item If the $\sigma$-algebras $\seq{\F_n}$ are independent, then $\F_{\infty}$ is $\PM$\-/trivial and $\QM$ is a.s.\ constant.
\end{enumerate}
\end{lemma}

\begin{proof}
Let $\P$ be the set of probability measures on $\X$ and fix $m\ge 1$. Note that $(\P,\BLM)$ is complete \citep[see, e.g.,][Theorem 1.12.4]{vanderVaart/Wellner:23}. For $k,l\ge m$, the kernels $\QM_k$ and $\QM_l$ are $\F_m^{\vee}$\-/measurable, so that $\dBL{\QM_k,\QM_l}$ is $\F_m^{\vee}$\-/measurable. Thus,
\[
	C\eqdef \{\omega: \seq{\QM_n(\omega,\csdot)} \text{ is $\BLM$\-/Cauchy}\} \in \F_m^{\vee}.
\]
Fix $\nu\in\P$ and for $n\ge 1$, define
\[
	\widetilde{\QM}_n(\omega,\csdot)\eqdef \ind_C(\omega)\QM_n(\omega,\csdot) + \ind_{C^c}(\omega)\nu(\csdot).
\]
The maps $\widetilde{\QM}_n$, $n\ge m$, are $\F_m^{\vee}$\-/measurable and converge pointwise to $\widetilde\QM$, defined as $\lim_n\QM_n$ on $C$ and $\nu$ on $C^{c}$. Therefore, $\widetilde\QM$ is $\F_m^{\vee}$\-/measurable. As $m$ is arbitrary and $\widetilde{\QM}$ does not depend on $m$, it is $\F_{\infty}$\-/measurable. Moreover, $\widetilde\QM$ is a probability kernel, and since $\QM_n\Rightarrow\QM$ a.s., we have $\PR{C}=1$ and $\widetilde\QM=\QM$ a.s. Parts (a) and (b) follow because $\F_m^{\vee}=\F^{\vee}$ and $\F_m^{\vee}=\F_m$ for every $m\ge1$, respectively, and (c) is Kolmogorov's zero--one law.
\end{proof}

\enlargethispage{0.5em}

A probability kernel is the conditional law of a random element given every $\F_n$ (on a suitable extension of the underlying probability space) if and only if it is $\F^{\wedge}$\-/measurable. Therefore, when $\F^{\wedge}$ is $\PM$\-/trivial, a fixed comparison is a.s.\ constant. Since a limit can be taken $\F_{\infty}$\-/measurable, it is the conditional law of a random element $X$ given some $\sigma$-algebra $\F\subset\F_{\infty}$.

We turn to $\FS$\-/weak convergence and summarize its properties, some of which are direct corollaries of Lemma \ref{lemma:merging}. Its characterization through conditional characteristic functions is provided in Remark \ref{rem:merging2}(2). The conditions of parts (a) and (b) are stated in terms of the limit kernel $\QM$. However, as noted after Lemma \ref{lemma:tail_join}, they can be interpreted as $X_n\to X$, $\FS$\-/weakly w.r.t.\ some $\F\subset\F_{\infty}$.

\begin{lemma}
\label{lemma:weak_conv}
Let $\seq{\QM_n}$ and $\QM$ be probability kernels on $\X$. Then,
\begin{enumerate}[label=\tn{(\alph*)}]
	\item \tn{(Converging together)} Let $\seq{X_n}$ and $\seq{Z_n}$ be sequences of $\X$\-/valued random elements. If $\QM_{X_n}^{\F_n}\Rightarrow\QM$ and $\varrho_{\X}(X_n,Z_n)\PRC{\FS}0$ a.s., then $\QM_{Z_n}^{\F_n}\Rightarrow\QM$.

	\item \tn{(Joint convergence)} Let $\seq{X_n}$ and $\seq{Y_n}$ be sequences of random elements of $\X$ and $\Y$, respectively. If $\QM_{X_n}^{\F_n}\Rightarrow\QM$ and $\varrho_{\Y}(Y_n,Y)\PRC{\FS}0$ a.s.\ for an $\F^{\wedge}$\-/measurable $Y$, then $\QM_{(X_n,Y_n)}^{\F_n}\Rightarrow\QM\otimes\delta_Y$ on $\X\times\Y$.

	\item \tn{(Continuous mapping)} Suppose $\QM_n\Rightarrow\QM$ and let $g:\X\to\Y$ be Borel and continuous outside a set $D_g$ with $\QM(\csdot,D_g)=0$ a.s. Then, $\QM_n\circ g^{-1}\Rightarrow \QM\circ g^{-1}$.

	\item \tn{(Cram\'er--Wold)} Let $\X=\R^d$ and let $D$ be a countable and dense subset of $\X$. Then $\QM_n\Rightarrow\QM$ if and only if for every $t\in D$, $\QM_n\circ(t^{\top}\csdot)^{-1}\Rightarrow\QM\circ(t^{\top}\csdot)^{-1}$.
\end{enumerate}
\end{lemma}

\begin{proof}
(a) This is Lemma \ref{lemma:merging}(a) with the constant comparison $\QM_n'\equiv\QM$.

\smallskip

\noindent(b) Since $Y$ is $\F_n$\-/measurable for every $n\ge 1$, $\QM_{(X_n,Y)}^{\F_n}=\QM_{X_n}^{\F_n}\otimes\delta_Y$ a.s., and since $\QM_{X_n}^{\F_n}\Rightarrow\QM$, the products converge weakly to $\QM\otimes\delta_Y$. Finally, (a) yields the result by noticing that the distance between $(X_n,Y_n)$ and $(X_n,Y)$ is $\varrho_{\Y}(Y_n,Y)\PRC{\FS}0$.

\smallskip

\noindent(c) Outside a single null set $\QM_n(\omega,\csdot)\Rightarrow\QM(\omega,\csdot)$ and $\QM(\omega,D_g)=0$, so the classical continuous mapping theorem applies at each such $\omega$.

\smallskip

\noindent(d) This is Lemma \ref{lemma:merging}(c) with the tightness hypothesis imposed on the constant comparison $\QM_n'\equiv\QM$ (as in Remark \ref{rem:merging2}(2)).
\end{proof}

\subsection{Merging and the Kolmogorov distance}
\label{sec:polya}

The following results relate merging of probability kernels with their closeness in the Kolmogorov distance.

\begin{lemma}
\label{lemma:seq_polya}
Let $\seq{\QM_n}$ and $\seq{\QM_n'}$ be probability kernels on $\R$ such that $\QM_n\mrg\QM_n'$ and $\seq{\QM_n'}$ is asymptotically equicontinuous. Then $\dK{\QM_n,\QM_n'}\to 0$ a.s. In particular, if $\seq{X_n}$ merges $\FS$\-/weakly with an $\FS$\-/asymptotically equicontinuous $\seq{Z_n}$, then $\dK{X_n,Z_n\mid\F_n}\to0$ a.s.
\end{lemma}

\begin{proof}
Let $F_n$ and $F_n'$ be the cdfs of $\QM_n$ and $\QM_n'$, respectively, and let $\beta_n\eqdef\dBL{\QM_n,\QM_n'}$. Fix $\delta>0$ and $\omega$ outside the union of the involved null sets. Using the asymptotic equicontinuity of $\seq{\QM_n'}$, choose $\epsilon>0$ with $\limsup_n\lC{\epsilon,\QM_n'(\omega)}\le\delta$. For $y\in\R$, let $f_y$ be the piecewise linear function defined as $1$ on $(-\infty,y]$, $0$ on $[y+\epsilon,\infty)$, and linear in between. Then $\norm{(\epsilon\wedge 1)f_y}_{\infty}\le 1$ and $\Lip((\epsilon\wedge 1)f_y)\le 1$, and hence
\[
	\abs{\int f_y\,d\QM_n(\omega,\csdot)-\int f_y\,d\QM_n'(\omega,\csdot)}\le (\epsilon\wedge 1)^{-1}\beta_n(\omega).
\]
Since for $y\in\R$, $\ind_{(-\infty,y]}\le f_{y}\le\ind_{(-\infty,y+\epsilon]}$, we get
\begin{align*}
	F_n(\omega,y)&\le F_n'(\omega,y+\epsilon)+(\epsilon\wedge 1)^{-1}\beta_n(\omega) \qtext{and} \\
	F_n(\omega,y)&\ge F_n'(\omega,y-\epsilon)-(\epsilon\wedge 1)^{-1}\beta_n(\omega)
\end{align*}
for every $y\in\R$. Therefore, $\abs{F_n(\omega,y)-F_n'(\omega,y)}\le \lC{\epsilon,\QM_n'(\omega)}+(\epsilon\wedge 1)^{-1}\beta_n(\omega)$. As $\epsilon$ was chosen independently of $y$, we have $\limsup_n\sup_{y\in\R}\abs{F_n(\omega,y)-F_n'(\omega,y)}\le\delta$. Sending $\delta\to 0$ implies the result.
\end{proof}

\begin{remark}
For weakly approaching sequences in $\R^d$ with a nonrandom comparison, the counterparts of Lemma \ref{lemma:seq_polya} are Lemma 9 and Corollary 1 of \cite{Belyaev:00}, where the comparison sequence is required to be tight, and the uniform continuity of its cdfs replaces asymptotic equicontinuity.
\end{remark}

\begin{lemma}
\label{lemma:tight_conv}
Let $\seq{\QM_n}$ and $\seq{\QM_n'}$ be probability kernels on $\R$. Suppose that $\dK{\QM_n,\QM_n'}\to 0 $ a.s.\ and $\{\QM_n'\}$ is asymptotically tight. Then, $\seq{\QM_n}$ is asymptotically tight and $\QM_n\mrg\QM_n'$. In particular, if $\dK{X_n,Z_n\mid\F_n}\to 0$ a.s.\ and $\seq{Z_n}$ is $\FS$\-/asymptotically tight, then $\seq{X_n}$ merges $\FS$\-/weakly with $\seq{Z_n}$.
\end{lemma}

\begin{proof}
Let $F_n$ and $F_n'$ be the cdfs of $\QM_n$ and $\QM_n'$, respectively, and let $\rho_n\eqdef \dK{\QM_n,\QM_n'}$. We show the asymptotic tightness of $\seq{\QM_n}$ first. For every $y>1$ and $\omega\in\Omega$,
\begin{align*}
	\QM_n(\omega,\{\abs{x}>y\})&\le 1-F_n(\omega,y)+F_n(\omega,-y) \\
	&\le 1-F_n'(\omega,y)+F_n'(\omega,-y)+2\rho_n(\omega) \\
	&\le\QM_n'(\omega,\{\abs{x}>y-1\})+2\rho_n(\omega).
\end{align*}
Since $\rho_n\to 0$ a.s., we have
\begin{align*}
	&\essinf_{y> 1}\limsup_{n\to\infty}\QM_n(\csdot,\{\absin{x}>y\}) \\
	&\qquad\le\essinf_{y> 0}\limsup_{n\to\infty}\QM_n'(\csdot,\{\absin{x}>y\})=0 \qtext{a.s.}
\end{align*}

As for merging, fix $\omega$ outside the union of the involved null sets. Let $\chi_y$ be the piecewise linear function defined as $1$ on $[-y,y]$, $0$ outside $[-y-1,y+1]$, and linear in between. For $f$ with $\norm{f}_{\infty}\le 1$ and $\Lip(f)\le 1$ define $h\eqdef f\chi_y$. Then $\norm{h}_{\infty}\le 1$, $\Lip(h)\le2$, and $h$ vanishes outside $[-y-1,y+1]$. Since $0\le 1-\chi_y\le\ind\{\abs{x}>y\}$ and $\abs{f}\le 1$,
\begin{align*}
	\abs{\int f\,d\QM_n(\omega,\csdot)-\int f\,d\QM_n'(\omega,\csdot)}&\le\abs{\int h\,d\QM_n(\omega,\csdot)-\int h\,d\QM_n'(\omega,\csdot)} \\
	&\quad+\QM_n(\omega,\{\abs{x}>y\})+\QM_n'(\omega,\{\abs{x}>y\}).
\end{align*}
Since $h$ is Lipschitz with compact support, we can write $h(x)=-\int h'(z)\ind\{z\ge x\}\,dz$,
so, by Fubini's theorem $\int h\,d\QM_n(\omega,\csdot)=-\int h'(z)F_n(\omega,z)\,dz$ and likewise for $\QM_n'$. Hence,
\begin{align*}
	\abs{\int h\,d\QM_n(\omega,\csdot)-\int h\,d\QM_n'(\omega,\csdot)}&=\abs{\int h'(z)\left(F_n(\omega,z)-F_n'(\omega,z)\right)\,dz} \\
	&\le\rho_n(\omega)\int\abs{h'(z)}\,dz\le 4(y+1)\rho_n(\omega).
\end{align*}
The bound is uniform over the bounded-Lipschitz ball, and $y$ was chosen independently of $f$. Therefore,
\begin{align*}
	\DBL{\QM_n(\omega,\csdot),\QM_n'(\omega,\csdot)}&\le\QM_n(\omega,\{\abs{x}>y\})+\QM_n'(\omega,\{\abs{x}>y\}) \\
	&\quad + 4(y+1)\rho_n(\omega).
\end{align*}
Since the probabilities on the right-hand side of the preceding display can be made arbitrarily small, we get $\dBL{\QM_n(\omega,\csdot),\QM_n'(\omega,\csdot)}\to 0$, implying the result.
\end{proof}

In the case of a fixed comparison with continuous cdf, the comparison is asymptotically tight and equicontinuous, so that the conditional form of P\'olya's theorem follows.

\begin{corollary}
\label{corr:weak_eq}
Let $\seq{\QM_n}$ and $\QM$ be probability kernels on $\R$ and suppose the cdf of $\QM$ is a.s.\ continuous. Then $\QM_n\Rightarrow\QM$ if and only if $\dK{\QM_n,\QM}\to 0$ a.s. In particular, if $\FD{X}{\F}$ is a.s.\ continuous, then \tn{(a)} $X_n\to X$, $\F$\-/weakly, if and only if $\dK{X_n,X\mid\F}\to 0$ a.s., and \tn{(b)} if $\F\subset\F^{\wedge}$ and for each $n\ge 1$, $\CI{X}{\F_n}{\F}$, then $X_n\to X$, $\FS$\-/weakly w.r.t.\ $\F$, if and only if $\dK{X_n,X\mid\F_n}\to 0$ a.s.
\end{corollary}

Lemmas \ref{lemma:seq_polya} and \ref{lemma:tight_conv} show how the equivalence from Corollary \ref{corr:weak_eq} adapts to a moving comparison sequence. Specifically, under asymptotic equicontinuity, merging implies Kolmogorov convergence, while under asymptotic tightness, convergence in the Kolmogorov distance implies merging. When both conditions hold jointly, the two notions coincide. Example \ref{ex:equicontinuity} below shows that asymptotic equicontinuity cannot be dropped from Lemma \ref{lemma:seq_polya}, while Example \ref{ex:tightness} demonstrates that asymptotic tightness cannot be dropped from Lemma \ref{lemma:tight_conv}.

\begin{example}
\label{ex:equicontinuity}
Let $\F_n\equiv\H_{\varnothing}$, $X_n\sim\nD{n^{-1}}{\sigma_n^2}$, and let $F_n'$ be the cdf of $\nD{0}{\sigma_n^2}$ with $\sigma_n=o(n^{-1})$. The variable $Z_n\eqdef X_n-n^{-1}$ has cdf $F_n'$, so that $\dBL{\QM_{X_n}^{\H_{\varnothing}},\QM_{Z_n}^{\H_{\varnothing}}}\le n^{-1}\to 0$ as $n\to\infty$. However, both $\dK{X_n,Z_n\mid \H_{\varnothing}}=2\Phi(1/(2n\sigma_n))-1$ and $\lC{\epsilon,Z_n\mid \H_{\varnothing}}=2\Phi(\epsilon/(2\sigma_n))-1$ converge to $1$ for any $\epsilon>0$. In particular, $\seq{Z_n}$ is not $\FS$\-/asymptotically equicontinuous.
\end{example}

\begin{example}
\label{ex:tightness}
Let $\QM_n'$ be the uniform mixture of the $n+1$ uniform laws on $I_k\eqdef[2k+1,2k+1+w_n]$, $0\le k\le n$, with $w_n\eqdef(n+1)^{-1}$. Define $\QM_n\eqdef\QM_n'\circ f^{-1}$, where $f(x)=x+1$. For a fixed $x>0$, at most $(x+1)/2$ of the intervals $I_k$ intersect $[-x,x]$, so that $\QM_n'([-x,x])\le (x+1)(2(n+1))^{-1}\to 0$ as $n\to\infty$, i.e., $\seq{\QM_n'}$ is not asymptotically tight. On the other hand, $\QM_n'$ has a density bounded by $1$, so that $\lC{\delta,\QM_n'}\le\delta$ for every $\delta>0$, and $\seq{\QM_n'}$ is asymptotically equicontinuous. Since $\QM_n((-\infty,x])=\QM_n'((-\infty,x-1])$ and an interval of length $1$ intersects at most one of the intervals $I_k$, we have $\dK{\QM_n,\QM_n'}=\sup_{x}\QM_n'((x-1,x])\le(n+1)^{-1}\to 0$ as $n\to\infty$. Consider the period-2 triangular wave $g$, which takes the value $1$ at odd integers and $0$ at even integers. It belongs to the unit bounded-Lipschitz ball. Averaging $g$ over each $I_k$ and each $f(I_k)$ yields $\int g\,d\QM_n'-\int g\,d\QM_n=1-w_n$, implying that $\dBL{\QM_n,\QM_n'}\ge 1-w_n\to 1$ as $n\to\infty$. \end{example}

\begin{remark}
(1) The difference $\DM_n\eqdef\QM_n-\QM_n'$ of two probability kernels is a signed kernel of total variation at most $2$, and both notions of closeness used here are norms of $\DM_n$. For finite signed measures $\mu$ on $\X$ and $\nu$ on $\R$, let $\norm{\mu}_{\mathrm{BL}}\eqdef\sup\{|\int f\,d\mu|:f\in \BL{\X}\}$ and $\norm{\nu}_{\mathrm{K}}\eqdef\sup_{x\in\R}\absin{\nu((-\infty,x])}$. (For example, $\QM_n\mrg\QM'_n$ if and only if $\norm{\DM_n}_{\mathrm{BL}}\to 0$ a.s.)

Neither norm dominates the other. In Example \ref{ex:equicontinuity} above, $\norm{\DM_n}_{\mathrm{BL}}\le n^{-1}$ while $\norm{\DM_n}_{\mathrm{K}}\to 1$, and in Example \ref{ex:tightness}, $\norm{\DM_n}_{\mathrm{K}}\le(n+1)^{-1}$ while $\liminf_n\norm{\DM_n}_{\mathrm{BL}}\ge 1$. Lemmas \ref{lemma:seq_polya} and \ref{lemma:tight_conv} provide conditions for transitioning between these norms. Converting $\norm{\csdot}_{\mathrm{K}}$ into $\norm{\csdot}_{\mathrm{BL}}$ relies on truncation, where asymptotic tightness ensures that the discarded tails are negligible. The reverse conversion requires smoothing an indicator over a small window, and the asymptotic equicontinuity of the comparison makes the smoothing error negligible. Specifically, for $\epsilon>0$, the proof of Lemma \ref{lemma:seq_polya} shows that
\begin{equation}
	\label{eq:dk_dbl_bound}
	\norm{\DM_n}_{\mathrm{K}}\le\lC{\epsilon,\QM_n'}+(\epsilon\wedge1)^{-1}\norm{\DM_n}_{\mathrm{BL}}.
\end{equation}
If the comparison has a density bounded by $c_n$, so that $\lC{\epsilon,\QM_n'}\le c_n\epsilon$, and $\norm{\DM_n}_{\mathrm{BL}}\le c_n$, choosing $\epsilon=(\norm{\DM_n}_{\mathrm{BL}}/c_n)^{1/2}$ yields $\norm{\DM_n}_{\mathrm{K}}\le 2(c_n\norm{\DM_n}_{\mathrm{BL}})^{1/2}$. Asymptotic equicontinuity is the anti-concentration hypothesis that replaces the bounded density when the comparison moves with $n$.

(2) Asymptotic equicontinuity can be checked through characteristic functions. By Esseen's concentration inequality, applied $\omega$-wise, $\lC{\epsilon,\QM}\le C\epsilon\int \ind\{\absin{t}\le1/\epsilon\}\absin{\widehat\QM(\csdot,t)}\,dt$ for an absolute constant $C$. Hence, if there is a jointly measurable $\psi:\Omega\times\R\to[0,\infty)$ with $\int\psi(\csdot,t)\,dt<\infty$ a.s.\ such that, for almost every $\omega\in\Omega$, $\absin{\widehat\QM_n(\omega,t)}\le\psi(\omega,t)$ for all $n$ and $t$, then $\lC{\epsilon,\QM_n}\le C\epsilon\int\psi(\csdot,t)\,dt$ a.s., and $\seq{\QM_n}$ is asymptotically equicontinuous. This holds, for example, in the case of Gaussian kernels with variances bounded away from zero.
\end{remark}

We conclude this section with a perturbation bound. Specifically, when a comparison sequence is asymptotically equicontinuous, perturbations that vanish in conditional probability do not affect the asymptotic behavior of the Kolmogorov distance.

\begin{lemma}
\label{lemma:nuisance}
Let $\seq{\G_n}$ be a sequence of sub-$\sigma$-algebras of $\H$ with $\F_n\subset\G_n$ for each $n\ge 1$, and let $\seq{X_n}$, $\seq{W_n}$, and $\seq{Z_n}$ be sequences of random variables. Then for any $\epsilon>0$,
\begin{align}
	\label{eq:nuisance_bound}
	\begin{aligned}
		\DK{X_n,Z_n\mid\G_n,\F_n}&\le\DK{W_n,Z_n\mid\G_n,\F_n} \\
		&\quad+\LC{\epsilon,Z_n\mid\F_n}+\PR{\abs{X_n-W_n}>\epsilon\mid\G_n} \qtext{a.s.}
	\end{aligned}
\end{align}
Consequently, if $\absin{X_n-W_n}\PRC{\FS}0$ and $\dK{W_n,Z_n\mid\G_n,\F_n}\PRC{\FS}0$ a.s., and $\seq{Z_n}$ is $\FS$\-/asymptotically equicontinuous, then $\dK{X_n,Z_n\mid\G_n,\F_n}\PRC{\FS}0$ a.s.
\end{lemma}

\begin{proof}
Fix $\epsilon>0$ and $r\in\R$. Since $\{X_n\le r\}\subset\{W_n\le r+\epsilon\}\cup\{\abs{X_n-W_n}>\epsilon\}$,
\begin{align*}
	\PR{X_n\le r\mid\G_n}&\le\PR{W_n\le r+\epsilon\mid\G_n}+\PR{\abs{X_n-W_n}>\epsilon\mid\G_n} \qtext{a.s.},
\shortintertext{while}
	\PR{W_n\le r+\epsilon\mid\G_n}&\le\FD{Z_n}{\F_n}(\csdot,r+\epsilon)+\DK{W_n,Z_n\mid\G_n,\F_n} \\
	&\le\FD{Z_n}{\F_n}(\csdot,r)+\lC{\epsilon,Z_n\mid\F_n}+\DK{W_n,Z_n\mid\G_n,\F_n} \qtext{a.s.}
\end{align*}
Thus, $\FD{X_n}{\G_n}(\csdot,r)-\FD{Z_n}{\F_n}(\csdot,r)$ is bounded by the right-hand side of \eqref{eq:nuisance_bound}. The symmetric argument applied to $\{W_n\le r-\epsilon\}\subset\{X_n\le r\}\cup\{\abs{X_n-W_n}>\epsilon\}$ bounds $\FD{Z_n}{\F_n}(\csdot,r)-\FD{X_n}{\G_n}(\csdot,r)$ by the same quantity. Taking the supremum over rational $r$ yields \eqref{eq:nuisance_bound}.

For the second claim, fix $\nu>0$. Then,
\begin{align*}
	\PR{\DK{X_n,Z_n\mid\G_n,\F_n}>3\nu\mid\F_n}&\le\PR{\DK{W_n,Z_n\mid\G_n,\F_n}>\nu\mid\F_n} \\
	&\quad+\ind\{\LC{\epsilon,Z_n\mid\F_n}>\nu\} \\
	&\quad+\nu^{-1}\PR{\abs{X_n-W_n}>\epsilon\mid\F_n} \qtext{a.s.}
\end{align*}
by the conditional Markov inequality and the tower property on the third term. Therefore,
\begin{align*}
	&\limsup_{n\to\infty}\PR{\DK{X_n,Z_n\mid\G_n,\F_n}>3\nu\mid\F_n} \\
	&\qquad\le \essinf_{\epsilon>0}\limsup_{n\to\infty}\ind\{\LC{\epsilon,Z_n\mid\F_n}>\nu\}=0 \qtext{a.s.} \qedhere
\end{align*}
\end{proof}

\section{Conditional Bootstrap}
\label{sec:cond_bootstrap}

\subsection{Setup}
\label{sec:setup}

In this section we present a self-contained formulation of the bootstrap under conditioning on a sequence of $\sigma$-algebras $\CS\eqdef\{\C_n\}$, representing the common shocks. Suppose that $\seq{X_n}$ is a sequence of random elements, where $X_n$ takes values in a Polish space $\X_n$, e.g., $\X_n=\X^n$ for a sample of size $n$ from $\X$. For a given $n\ge1$, let $\theta_n$ be a $\C_n$\-/measurable parameter and let
\[
	T_n(\theta_n)\eqdef\T_n(X_n,\theta_n),
\]
where $\T_n$ is a measurable real-valued function, denote a statistic used to conduct inference on $\theta_n$ based on a realization of $X_n$. We denote the conditional cdf of $T_n$ given $\C_n$ by $\FD{n}{}$. We provide sufficient conditions for the conditional first-order consistency of resampling estimators of $\FD{n}{}$. Specifically, let $\G_n\eqdef\C_n\vee\sigma(X_n)$ and let $X_n^{*}$ be a pseudo-sample drawn using a realization of $X_n$. Then the bootstrap counterpart of $T_n$ is $T_n^{*}\eqdef\T_m(X_n^{*},\theta_n^{*})$, where $m\equiv m(n)$ is the size of $X_n^{*}$ and $\theta_n^{*}$ is an estimator of $\theta_n$ computed from $X_n$. The conditional cdf $\FD{n}{*}$ of $T_n^{*}$ given $\G_n$ is used to approximate $\FD{n}{}$. By enlarging $(\Omega,\H,\PM)$ if necessary, we assume that the resampling weights that generate $X_n^{*}$ from $X_n$ and the comparison variables of Section \ref{sec:consistency} are defined on this probability space as well.

In applications, the statistic frequently involves a further unknown $\C_n$\-/measurable \emph{nuisance} $\vartheta_n$, e.g., a variance or scale. Its presence does not affect the results of Section \ref{sec:consistency} unless $\vartheta_n$ needs to be replaced by an estimator. We write $\widetilde{T}_n$ for the \emph{feasible} statistic obtained from $T_n$ by substituting an estimator $\widehat\vartheta_n$ for $\vartheta_n$, and $\widetilde{T}_n^{*}$ for its bootstrap counterpart, in which the substituted value $\widehat\vartheta_n^{\,*}$ may be $\widehat\vartheta_n$ itself, held fixed across draws, or a value recomputed from each pseudo-sample.

\subsection{Consistency}
\label{sec:consistency}

In the unconditional case (i.e., $\C_n\equiv\H_{\varnothing}$), bootstrap consistency is typically established by showing that $T_n$ and $T_n^{*}$ are close in the Kolmogorov distance to a common comparison, usually a weak limit of $T_n$. If this limit has a continuous cdf, the conditional quantiles of $\FD{n}{*}$ approximate those of $\FD{n}{}$, i.e., $\FD{n}{}(c_n^{*}(\alpha))\to\alpha$ in probability, where $c_n^{*}(\alpha)\eqdef\FD{n}{*}{}^{-}(\csdot,\alpha)$ is the conditional $\alpha$\-/quantile of $\FD{n}{*}$. This typically happens when the statistic $T_n$ is pivotal. However, the cdfs of a non-pivotal statistic, e.g., one whose conditional law retains dependence on a feature of the data-generating process that is hard to estimate consistently, need not converge,\footnote{
	Consider, e.g., the linearized statistic $T_n'$ in \eqref{eq:linear_stat}, whose law depends on $\theta_n$. It need not have a nondegenerate weak limit when $\seq{\theta_n}$ does not converge. Section \ref{sec:example_fd} provides an example in which the conditional laws of $T_n$ do not converge even after $\C_n$\-/measurable affine transformations.
}
and the comparison must vary with $n$. Without a continuous limit, closeness of $T_n^{*}$ and $T_n$ in the Kolmogorov distance does not necessarily yield $\FD{n}{}(c_n^{*}(\alpha))\to\alpha$. However, as shown below, it suffices for $T_n$ and $T_n^{*}$ to be close to a common sequence of comparisons whose conditional laws given $\C_n$ have asymptotically vanishing atoms.

Throughout this subsection, we fix a sequence $\seq{S_n}$ of comparison random variables and define
\[
	\Delta_n\eqdef\DK{T_n,S_n\mid\C_n}
	\qtext{and}\quad
	\Delta_n^{*}\eqdef\DK{T_n^{*},S_n\mid\G_n,\C_n}.
\]
These quantities measure how well the conditional law of $S_n$ given $\C_n$ approximates the conditional laws of $T_n$ and $T_n^{*}$ given $\C_n$ and $\G_n$, respectively.

\begin{lemma}
\label{lemma:coverage_bound}
For every $n\ge 1$, $\alpha\in(0,1)$, and $\eta\in(0,1)$,
\begin{align}
	\label{eq:coverage_bound}
	\begin{aligned}
		&\abs{\PR{T_n\le c_n^{*}(\alpha)\mid\C_n}-\alpha} \\
		&\qquad\le\eta+\lCZ{S_n\mid\C_n}+\Delta_n+\PR{\Delta_n^{*}>\eta\mid\C_n} \qtext{a.s.}
	\end{aligned}
\end{align}
\end{lemma}

\begin{proof}
Let $a_n\eqdef\lCZ{S_n\mid \C_n}$ and let $F_n'\eqdef \FD{S_n}{\C_n}$, whose conditional $\beta$-quantile is denoted by $c_n'(\beta)$. We use the following properties of generalized inverses \citep[see, e.g.,][Proposition 1]{Embrechts/Hofert:13}: for a cdf $F$ and $\beta\in(0,1)$, (P1) $F(F^{-}(\beta)-)\le\beta\le F(F^{-}(\beta))$, and (P2) $F(x)\ge\beta$ if and only if $x\ge F^{-}(\beta)$. Since the jumps of $F_n'$ are a.s.\ bounded by $a_n$, (P1) yields
\begin{equation}
	\label{eq:quantile_props}
	\beta\le F_n'(\csdot,c_n'(\beta))\le\beta+a_n \qtext{and}\quad \FD{n}{*}(\csdot,c_n^{*}(\beta))\ge\beta \qtext{a.s.}
\end{equation}


In addition, for a random variable $Z$, $\PR{T_n\le Z\mid\C_n}=\FD{n}{}(\csdot,Z)$ a.s., when $Z$ is $\C_n$\-/measurable \citep[see, e.g.,][Theorem 8.5]{Kallenberg:21:FMP}, and
\begin{align}
	\label{eq:random_dk}
	\begin{aligned}
		&\absin{\FD{n}{}(\csdot,Z)-F_n'(\csdot,Z)}\le\Delta_n \qtext{and} \\
		&\absin{\FD{n}{*}(\csdot,Z)-F_n'(\csdot,Z)}\le\Delta_n^{*} \qtext{a.s.}
	\end{aligned}
\end{align}

Fix $\alpha\in(0,1)$ and $\eta\in(0,1)$, and let $\Xi_n\eqdef\{\Delta_n^{*}\le\eta\}$. For the upper bound, suppose that $\alpha+\eta<1$. By \eqref{eq:quantile_props} and \eqref{eq:random_dk},
\[
	\FD{n}{*}(\csdot,c_n'(\alpha+\eta))\ge F_n'(\csdot,c_n'(\alpha+\eta))-\Delta_n^{*}\ge\alpha \qtext{a.s.\ on }\Xi_n,
\]
which implies that $c_n^{*}(\alpha)\le c_n'(\alpha+\eta)$ a.s.\ on $\Xi_n$ by (P2). Hence,
\begin{align*}
	\PR{T_n\le c_n^{*}(\alpha)\mid\C_n}&\le\FD{n}{}(\csdot,c_n'(\alpha+\eta))+\PR{\Delta_n^{*}>\eta\mid\C_n} \\
	&\le F_n'(\csdot,c_n'(\alpha+\eta))+\Delta_n+\PR{\Delta_n^{*}>\eta\mid\C_n} \\
	&\le\alpha+\eta+a_n+\Delta_n+\PR{\Delta_n^{*}>\eta\mid\C_n} \qtext{a.s.}
\end{align*}
The last line is at least $1$ a.s.\ if $\alpha+\eta\ge 1$, i.e., the bound is trivial. For the lower bound, suppose that $\alpha>\eta$. Similarly,
\[
	F_n'(\csdot,c_n^{*}(\alpha))\ge\FD{n}{*}(\csdot,c_n^{*}(\alpha))-\Delta_n^{*}\ge\alpha-\eta \qtext{a.s.\ on }\Xi_n,
\]
i.e., $c_n^{*}(\alpha)\ge c_n'(\alpha-\eta)$ a.s.\ on $\Xi_n$. Thus,
\begin{align*}
	\PR{T_n\le c_n^{*}(\alpha)\mid\C_n}&\ge\FD{n}{}(\csdot,c_n'(\alpha-\eta))-\PR{\Delta_n^{*}>\eta\mid\C_n} \\
	&\ge F_n'(\csdot,c_n'(\alpha-\eta))-\Delta_n-\PR{\Delta_n^{*}>\eta\mid\C_n} \\
	&\ge\alpha-\eta-\Delta_n-\PR{\Delta_n^{*}>\eta\mid\C_n} \qtext{a.s.}
\end{align*}
If $\alpha\le\eta$, this bound holds trivially. Combining these bounds yields \eqref{eq:coverage_bound}.
\end{proof}

Combining Lemma \ref{lemma:coverage_bound} with the Ky Fan distance of $\Delta_n^{*}$ from zero under its conditional law given $\C_n$, defined by
\[
	\kappa_n\eqdef\inf\left\{\epsilon>0:\FD{\Delta_n^{*}}{\C_n}(\csdot,\epsilon)\ge1-\epsilon\right\},
\]
we get a finite-sample uniform bound on the coverage error in terms of the two Kolmogorov distances, such as those in the form of Berry--Esseen bounds.\footnote{
	Since $\epsilon\mapsto\FD{\Delta_n^{*}}{\C_n}(\omega,\epsilon)$ is nondecreasing and right-continuous, $\{\kappa_n\le t\}=\big\{\FD{\Delta_n^{*}}{\C_n}(\csdot,t)\ge 1-t\big\}$ for $t>0$.
	Therefore, $\kappa_n$ is $\C_n$\-/measurable.
}

\begin{theorem}
\label{thm:coverage_bound}
For every $n\ge 1$,
\[
	\esssup_{\alpha\in(0,1)}\abs{\PR{T_n\le c_n^{*}(\alpha)\mid\C_n}-\alpha}\le\lCZ{S_n\mid\C_n}+\Delta_n+2\kappa_n \qtext{a.s.}
\]
\end{theorem}

\begin{proof}
Fix $\alpha\in(0,1)$. Replacing $\PR{\Delta_n^{*}>\eta\mid\C_n}$ with $1-\FD{\Delta_n^{*}}{\C_n}(\csdot,\eta)$ in \eqref{eq:coverage_bound}, we see that the bound holds for all $\eta\in [0,1)$ a.s.\ by the right-continuity of the regular conditional cdfs. Thus, it holds a.s.\ for a $\C_n$\-/measurable random variable $\eta\in [0,1)$ a.s., in particular for $\eta=\kappa_n\ind\{\kappa_n<1/2\}$. Since $1-\FD{\Delta_n^{*}}{\C_n}(\csdot,\kappa_n)\le\kappa_n$, with this choice of $\eta$, we get
\begin{equation}
	\label{eq:coverage_bound2}
	\abs{\PR{T_n\le c_n^{*}(\alpha)\mid\C_n}-\alpha}
	\le \lCZ{S_n\mid\C_n}+\Delta_n+2\kappa_n
\end{equation}
a.s.\ on $\{0\le \kappa_n<1/2\}$. The result follows because \eqref{eq:coverage_bound2} does not depend on $\alpha$, and it is trivial on $\{\kappa_n\ge 1/2\}$.
\end{proof}

\begin{remark}
\label{rem:coverage}
When $\lCZ{T_n\mid\C_n}\to0$ a.s., one may take $S_n=T_n$. In this case, $\Delta_n=0$, $\Delta_n^{*}=\DK{T_n^{*},T_n\mid\G_n,\C_n}$, and the bound of Theorem \ref{thm:coverage_bound} becomes $\lCZ{T_n\mid\C_n}+2\kappa_n$. A comparison with a known conditional law, e.g., a Gaussian one, is useful because $\Delta_n$ and $\Delta_n^{*}$ can be bounded using Berry--Esseen-type results.
\end{remark}

Theorem \ref{thm:coverage_bound} converts explicit bounds on $\Delta_n$ and on the conditional distribution of $\Delta_n^{*}$ into a bound on the conditional coverage error. Convergence $\Delta_n^{*}\PRC{\CS}0$ a.s.\ is equivalent to $\kappa_n\to0$ a.s. Indeed, if $\PR{\Delta_n^{*}>\epsilon\mid\C_n}\to 0$ a.s.\ for each $\epsilon>0$, then $\kappa_n\le\epsilon$ eventually a.s., and, conversely, $\PR{\Delta_n^{*}>\epsilon\mid\C_n}\le\kappa_n$ a.s.\ on $\{\kappa_n<\epsilon\}$. Moreover, if $r_n$ and $\beta_n$ are nonnegative $\C_n$\-/measurable random variables such that $\PR{\Delta_n^{*}>r_n\mid\C_n}\le\beta_n$ a.s., then $\kappa_n\le r_n\vee\beta_n$ a.s. In particular, by the conditional Markov inequality, for every $p>0$,
\[
	\kappa_n\le\left(\E[(\Delta_n^{*})^p\mid\C_n]\right)^{1/(p+1)} \qtext{a.s.}
\]

\begin{theorem}
\label{thm:bootstrap_consistency1}
Suppose that $\lCZ{S_n\mid\C_n}\to 0$ a.s.\ and the following conditions hold:
\begin{enumerate}[label=\tn{(\alph*)}]
	\item $\Delta_n\to 0$ a.s.

	\item $\Delta_n^{*}\PRC{\CS}0$ a.s.
\end{enumerate}
Then $\dK{T_n^{*},T_n\mid\G_n,\C_n}\PRC{\CS}0$ a.s.\ and
\[
	\esssup_{\alpha\in(0,1)}\abs{\PR{T_n\le c_n^{*}(\alpha)\mid\C_n}-\alpha}\to0 \qtext{a.s.}
\]
\end{theorem}

\begin{proof}
First, by the triangle inequality, $\DK{T_n^{*},T_n\mid\G_n,\C_n}\le\Delta_n+\Delta_n^{*}$ a.s., and, since $\Delta_n$ is $\C_n$\-/measurable, the right-hand side vanishes in conditional probability by (a) and (b) as $n\to\infty$. The second assertion follows from Theorem \ref{thm:coverage_bound} because $\kappa_n\to 0$ a.s.\ by (b).
\end{proof}

\begin{remark}
(1) A variant of this result can be found in \cite{Chernozhukov/Chetverikov/Kato:13}. The comparison there is the Gaussian analog of a maximum that need not have a limit distribution, and their argument relies on an anti-concentration bound for maxima of Gaussian vectors, which is a bound on the L\'evy concentration function of the comparison. Section \ref{sec:example_hd} illustrates this setting conditionally on a common shock.

(2) The bootstrap critical values $\{c_n^{*}(\alpha):\alpha\in(0,1)\}$ yield asymptotically correct unconditional coverage. In particular, by the dominated convergence theorem,
\begin{align*}
	&\sup_{\alpha\in(0,1)}\abs{\PR{T_n\le c_n^{*}(\alpha)}-\alpha}\\
	&\qquad\le\E\esssup_{\alpha\in(0,1)}\abs{\PR{T_n\le c_n^{*}(\alpha)\mid\C_n}-\alpha}\to 0
\end{align*}
as $n\to\infty$.
\end{remark}

\begin{definition}
\label{def:bootstrap_consistency}
We say that $\FD{n}{*}$ is conditionally $\KM$\-/consistent given $\CS$ if the conclusion of Theorem \ref{thm:bootstrap_consistency1} holds.
\end{definition}

Condition (b) of Theorem \ref{thm:bootstrap_consistency1} is typically easy to verify, for instance when the elements of $X_n^{*}$ are conditionally i.i.d.\ given $\G_n$, while establishing condition (a) can be demanding. However, when $\seq{T_n}$ and $\seq{S_n}$ merge $\CS$-weakly and the comparison is $\CS$\-/asymptotically equicontinuous, the bound in \eqref{eq:dk_dbl_bound} leads to the following corollary.

\begin{corollary}
\label{corr:bootstrap_merging}
Suppose that $\seq{S_n}$ is $\CS$\-/asymptotically equicontinuous and the following conditions hold:
\begin{enumerate}[label=\tn{(\alph*)}]
	\item $\seq{T_n}$ merges $\CS$\-/weakly with $\seq{S_n}$.

	\item $\dBL{\QM_{T_n^{*}}^{\G_n},\QM_{S_n}^{\C_n}}\PRC{\CS}0$ a.s.
\end{enumerate}
Then $\FD{n}{*}$ is conditionally $\KM$\-/consistent given $\CS$.
\end{corollary}

\begin{proof}
First, the asymptotic equicontinuity of $\seq{S_n}$ implies that $\lCZ{S_n\mid\C_n}\to0$ a.s. It remains to check conditions (a) and (b) of Theorem \ref{thm:bootstrap_consistency1}. Its condition (a) follows from Lemma \ref{lemma:seq_polya} applied to $\QM_{T_n}^{\C_n}$ and $\QM_{S_n}^{\C_n}$. Let $\beta_n\eqdef\dBL{\QM_{T_n^{*}}^{\G_n},\QM_{S_n}^{\C_n}}$ and fix $\nu>0$ and $\epsilon>0$. Using \eqref{eq:dk_dbl_bound} applied to the same kernels,
\[
	\PR{\Delta_n^{*}>2\nu\mid\C_n} \le\ind\{\LC{\epsilon,S_n\mid\C_n}>\nu\}+\PR{\beta_n>\nu(\epsilon\wedge1)\mid\C_n} \qtext{a.s.}
\]
Therefore, (b) implies that
\begin{align*}
	&\limsup_{n\to\infty}\PR{\Delta_n^{*}>2\nu\mid\C_n} \\
	&\qquad\le \essinf_{\epsilon>0}\limsup_{n\to\infty}\ind\{\LC{\epsilon,S_n\mid\C_n}>\nu\}=0 \qtext{a.s.},
\end{align*}
i.e., condition (b) of Theorem \ref{thm:bootstrap_consistency1} holds.
\end{proof}

Condition (b) of Theorem \ref{thm:bootstrap_consistency1} and Corollary \ref{corr:bootstrap_merging} is stated in conditional probability given $\C_n$ because the coverage in Theorem \ref{thm:bootstrap_consistency1} requires averaging over the data conditionally on the shock. The arguments of Section \ref{sec:cond_conv}, however, typically provide almost sure convergence of the bootstrap distances, e.g., $\Delta_n^{*}\to 0$ a.s. The next lemma provides conditions under which almost sure convergence implies convergence in conditional probability given $\CS$.

\begin{lemma}
\label{lemma:as_to_cond}
Let $\seq{\beta_n}$ be a sequence of random variables such that $\beta_n\to0$ a.s. Suppose that either \tn{(i)} $\CS$ is nondecreasing, or \tn{(ii)} for each $n\ge 1$, $\beta_n$ is $\G_n$\-/measurable and $\CI{\G_n}{\C^{\vee}}{\C_n}$. Then $\beta_n\PRC{\CS}0$ a.s. In particular, condition \tn{(b)} of Theorem \ref{thm:bootstrap_consistency1} holds if $\Delta_n^{*}\to 0$ a.s., and condition \tn{(b)} of Corollary \ref{corr:bootstrap_merging} holds if $\QM_{T_n^{*}}^{\G_n}\mrg\QM_{S_n}^{\C_n}$.
\end{lemma}

\begin{proof}
Fix $\epsilon>0$. Since $\beta_n\to 0$ a.s., $\ind\{\absin{\beta_n}>\epsilon\}\to 0$ a.s. Under (i), $\PR{\absin{\beta_n}>\epsilon\mid\C_n}\to 0$ a.s.\ by Hunt's lemma \citep[see, e.g.,][Theorem V.45]{Dellacherie/Meyer:82:ProbabilitiesB}. Under (ii), since $\C_n\subset\C^{\vee}$, $\PR{A\mid\C_n}=\PR{A\mid\C^{\vee}}$ a.s.\ for each $A\in\G_n$. Thus, by the conditional dominated convergence theorem we get
\[
	\PR{\absin{\beta_n}>\epsilon\mid\C_n}=\PR{\absin{\beta_n}>\epsilon\mid\C^{\vee}}\to 0 \qtext{a.s.}
\]
The last assertion follows because $\C_n\subset\G_n$, and thus both distances are $\G_n$\-/measurable.
\end{proof}

When the comparison law has a fixed, continuous conditional cdf given a sub-$\sigma$-algebra of $\C^{\wedge}$, we obtain the following case.

\begin{corollary}
\label{corr:bootstrap_consistency1}
Let $\C\subset\C^{\wedge}$ and suppose that $S_n\equiv S$ is such that $\FD{S}{\C}$ is a.s.\ continuous and $\CI{S}{\C_n}{\C}$ for each $n\ge1$. Furthermore, assume that the following conditions hold:
\begin{enumerate}[label=\tn{(\alph*)}]
	\item $T_n\to S$, $\CS$\-/weakly w.r.t.\ $\C$.

	\item $\Delta_n^{*}\PRC{\CS}0$ a.s.
\end{enumerate}
Then $\FD{n}{*}$ is conditionally $\KM$\-/consistent given $\CS$.
\end{corollary}

\begin{proof}
Since $\CI{S}{\C_n}{\C}$ and $\C\subset\C_n$, $\QM_S^{\C_n}=\QM_S^{\C}$ a.s., so that $\lCZ{S\mid\C_n}=0$ a.s., and Corollary \ref{corr:weak_eq}(b) yields $\Delta_n\to0$ a.s. The result follows from Theorem \ref{thm:bootstrap_consistency1}.
\end{proof}

Corollary \ref{corr:bootstrap_consistency1} requires the conditional law of the comparison given $\C_n$ to be the same for every $n$, i.e., $\C^{\wedge}$\-/measurable, which by the discussion following Lemma \ref{lemma:tail_join} forces an a.s.\ constant law whenever $\C^{\wedge}$ is $\PM$\-/trivial. When the conditioning fields are nondecreasing, any limit can be taken $\C^{\vee}$\-/measurable by Lemma \ref{lemma:tail_join}(a). As shown in the following result, a genuinely random comparison becomes available as soon as $\C^{\vee}$ is not $\PM$\-/trivial.

\begin{corollary}
\label{corr:bootstrap_consistency2}
Suppose that $\CS$ is nondecreasing and $S_n\equiv S$ is a constant sequence such that $\FD{S}{\C_{\infty}}$ is a.s.\ continuous. Then, $\dK{\QM_S^{\C_n},\QM_S^{\C_{\infty}}}\to 0$ a.s., and the following are equivalent:
\begin{enumerate}[label=\tn{(\roman*)}]
	\item $\Delta_n\to 0$ a.s.\ and $\Delta_n^{*}\PRC{\CS}0$ a.s.

	\item $\dK{\QM_{T_n}^{\C_n},\QM_S^{\C_{\infty}}}\to 0$ a.s.\ and $\dK{\QM_{T_n^{*}}^{\G_n},\QM_S^{\C_{\infty}}}\PRC{\CS}0$ a.s.
\end{enumerate}
Under either \tn{(i)} or \tn{(ii)}, $\FD{n}{*}$ is conditionally $\KM$\-/consistent given $\CS$.
\end{corollary}

\begin{proof}
Let $\rho_n\eqdef\dK{\QM_S^{\C_n},\QM_S^{\C_{\infty}}}$. Since $\C_n\subset\C_{\infty}$ and $[x,x+\epsilon]\subset[q,q+2\epsilon]$ for any $x\in\R$, $\epsilon>0$, and any rational $q\in(x-\epsilon,x]$,
\begin{align*}
	\LC{\epsilon,\QM_{S}^{\C_n}}\le\sup_{q\in\Q}\QM_S^{\C_n}(\csdot,[q,q+2\epsilon])\le\E\left[\LC{2\epsilon,\QM_{S}^{\C_{\infty}}}\mid \C_n\right] \qtext{a.s.}
\end{align*}
Since $\FD{S}{\C_{\infty}}$ is a.s.\ continuous, sending $\epsilon$ along a sequence converging to zero, the right-hand side vanishes a.s.\ by the conditional dominated convergence theorem. Thus, each $\FD{S}{\C_n}$ is a.s.\ continuous because jumps of $\FD{S}{\C_n}$ are bounded by $\lC{\epsilon,\QM_{S}^{\C_n}}$ for any $\epsilon>0$.

In addition, $\FD{S}{\C_n}(\csdot,q)\to\FD{S}{\C_{\infty}}(\csdot,q)$ a.s.\ for each $q\in\Q$ by L\'evy's upward theorem, applied to the constant sequence $\ind\{S\le q\}$. Since the limit is a.s.\ continuous, the classical P\'olya argument extends convergence on a dense set to uniform convergence, implying $\rho_n\to 0$ a.s. Therefore, $\rho_n\PRC{\CS}0$ a.s.\ by Lemma \ref{lemma:as_to_cond}(i). The equivalence of (i) and (ii) follows from the triangle inequality.

For the final assertion, note that $\lCZ{S\mid\C_n}\to 0$ a.s.\ because $\lCZ{S\mid\C_n}=0$ a.s.\ for each $n\ge 1$. Thus, the result follows from Theorem \ref{thm:bootstrap_consistency1}, applied with $S_n\equiv S$.
\end{proof}

\section{Smooth Statistics}
\label{sec:sfm}

In this section we consider the case when the statistic depends smoothly on the quantities that have to be estimated. The results of this part rest on the following observation. Suppose that a tractable auxiliary statistic satisfies the conditions of Theorem \ref{thm:bootstrap_consistency1} against an asymptotically equicontinuous comparison sequence, and suppose that the difference between the auxiliary statistic and the statistic of interest vanishes in conditional probability. Then, by Lemma \ref{lemma:nuisance}, the required conditions hold for the statistic of interest as well. We assume that $\X$ and $\Y$ are separable Banach spaces. Abusing notation, we write $\normin{\csdot}$ for their norms and for the operator norms of bounded linear and bilinear maps defined on them. The spaces and maps below may also depend on $n$ (see Remark \ref{rem:sfm_varying}). To have all the quantities in this section defined for every $\omega\in\Omega$, we set these maps and their derivatives to zero outside the sets on which they are defined.

\subsection{Smoothness in the estimator}
\label{sec:linearization}

Suppose that each $\theta_n\in\Theta$ a.s.\ for some Borel set $\Theta\subset\X$. We first consider the case in which the statistic $T_n$ takes the form
\[
	T_n(\theta_n)=\ell\big(\tau_n(\phi(\widehat{\theta}_n)-\phi(\theta_n))\big),
\]
where $\phi:\Theta\to\Y$ is continuous, $\widehat{\theta}_n\in \Theta$ a.s.\ is a consistent estimator of $\theta_n$ (in the sense of Definition \ref{def:CCP}, i.e., $\normin{\widehat{\theta}_n-\theta_n}\PRC{\CS}0$ a.s.), and $\tau_n>0$ a.s.\ is a $\C_n$\-/measurable normalizing coefficient. The function $\ell:\Y\to\R$ is Lipschitz and yields the real statistic required by Section \ref{sec:cond_bootstrap}. In particular, the \textit{smooth function model} (see, e.g., \citeay{Hall:92:Bootstrap}, Section 2.4 and \citeay{Lahiri:03:Resampling}, Section 4.2) falls into this case. We assume that for almost every $\omega\in\Omega$, $\phi$ is Fr\'echet differentiable at $\theta_n(\omega)$ relative to $\Theta$, i.e., there is a bounded linear map $D\phi(\theta_n(\omega)):\X\to\Y$ such that
\[
	\norm{\phi(\theta)-\phi(\theta_n(\omega))-D\phi(\theta_n(\omega))(\theta-\theta_n(\omega))}=o(\norm{\theta-\theta_n(\omega)})
\]
as $\theta\to\theta_n(\omega)$ in $\Theta$. When $\theta_n(\omega)$ is not an interior point of $\Theta$, such a derivative need not be unique, and $D\phi(\theta_n(\omega))$ denotes a fixed choice. By the convention above, $D\phi(\theta_n(\omega))=0$ for the remaining $\omega\in\Omega$, and we assume that $\omega\mapsto D\phi(\theta_n(\omega))x$ is $\C_n$\-/measurable for every $x\in\X$.\footnote{
	This holds if $\phi$ is the restriction of a Fr\'echet differentiable map on an open $\Theta'\supset\Theta$. In this case, $D\phi(\theta)x=\lim_{m\to\infty}m(\phi(\theta+x/m)-\phi(\theta))$ is a Borel function of $\theta\in\Theta'$ for every $x\in\X$, and $D\phi(\theta_n)=0$ on $\{\theta_n\notin\Theta'\}$ by the same convention.
}
The resampling version of $T_n$ is
\[
	T_n^{*}=\ell\big(\tau_n^{*}(\phi(\widehat{\theta}_n^{*})-\phi(\theta_n^{*}))\big),
\]
where $\theta_n^{*}\in\Theta$ a.s.\ is a consistent estimator of $\theta_n$, which may differ from $\widehat{\theta}_n$, $\widehat{\theta}_n^{*}\in\Theta$ a.s.\ is the bootstrap analog of $\widehat{\theta}_n$, and $\tau_n^{*}>0$ a.s.\ is the $\G_n$\-/measurable bootstrap counterpart of $\tau_n$, which may depend on the data, as with a data-driven choice of $m$ in the $m$-out-of-$n$ bootstrap of \cite{Bickel/Sakov:08}.

Let $\xi_n\eqdef\tau_n(\widehat{\theta}_n-\theta_n)$, $\xi_n^{*}\eqdef\tau_n^{*}(\widehat{\theta}_n^{*}-\theta_n^{*})$, and $\rho_n\eqdef\tau_n\wedge\tau_n^{*}$. For $t\ge 0$, define
\[
	\varpi_n(t)\eqdef\sup\left\{\frac{\norm{\phi(\theta)-\phi(\theta_n)-D\phi(\theta_n)(\theta-\theta_n)}}{\norm{\theta-\theta_n}}:\theta\in\Theta,\ 0<\norm{\theta-\theta_n}<t\right\},
\]
the first-order remainder modulus of $\phi$ at $\theta_n$, with $\sup{\varnothing}\eqdef 0$.\footnote{
	Let $\Theta_0$ be a countable dense subset of $\Theta$. Since $\phi$ is continuous and the set $\{\theta\in\Theta:0<\normin{\theta-\theta_n}<t\}$ is relatively open in $\Theta$, the supremum may be taken over $\theta\in\Theta_0$. Hence, $\varpi_n(R)$ is $\F$\-/measurable for every $\sigma$-algebra $\F$ with respect to which $\theta_n$, $\omega\mapsto D\phi(\theta_n(\omega))x$ for $x\in\X$, and $R\ge0$ are measurable.
}
We impose the following assumption on these quantities.

\begin{assumption}
\label{assum:smooth_tightness}
The sequences $\seq{\normin{\xi_n}}$, $\seq{\normin{\xi_n^{*}}}$, and $\seq{\tau_n^{*}\normin{\theta_n^{*}-\theta_n}}$ are $\CS$\-/asymptotically tight, and $\varpi_n(M/\rho_n)\PRC{\CS}0$ a.s.\ for each $M\in [0,\infty)$.
\end{assumption}

Consider the linearized statistics
\begin{equation}
\label{eq:linear_stat}
	T_n'\eqdef \ell\big(D\phi(\theta_n)\xi_n\big) \qtext{and}\quad T_n'^{*}\eqdef \ell\big(D\phi(\theta_n)\xi_n^{*}\big),
\end{equation}
and let $V_n\eqdef\{(\theta_n,\widehat{\theta}_n, \theta_n^{*},\widehat{\theta}_n^{*})\in\Theta^4\}\cap\{(\tau_n,\tau_n^{*})\in (0,\infty)^2\}$. The following result shows that it suffices to find a suitable approximation of these statistics in order to apply Theorem \ref{thm:bootstrap_consistency1} to this setup.

\begin{lemma}
\label{lemma:bootstrap_consistency2}
Suppose that Assumption \ref{assum:smooth_tightness} holds and let $\seq{S_n'}$ be a $\CS$\-/asymptotically equicontinuous sequence of random variables satisfying
\begin{enumerate}[label=\tn{(\alph*)}]
	\item $\dK{T_n',S_n'\mid\C_n}\to 0$ a.s., and

	\item $\dK{T_n'^{*},S_n'\mid\G_n,\C_n}\PRC{\CS}0$ a.s.
\end{enumerate}
Then $\dK{T_n,S_n'\mid\C_n}\to 0$ and $\dK{T_n^{*},S_n'\mid\G_n,\C_n}\PRC{\CS}0$ a.s.
\end{lemma}

\begin{proof}
Fix $M\in(0,\infty)$, let $t_n\eqdef3M/\rho_n$, and for $\theta\in\Theta$, let $r_n(\theta)\eqdef\phi(\theta)-\phi(\theta_n)-D\phi(\theta_n)(\theta-\theta_n)$, so that $\normin{r_n(\theta)}\le\varpi_n(t_n)\normin{\theta-\theta_n}$ whenever $\normin{\theta-\theta_n}<t_n$. Then
\[
	\tau_n(\phi(\widehat\theta_n)-\phi(\theta_n))-D\phi(\theta_n)\xi_n=\tau_nr_n(\widehat\theta_n),
\]
and, by the linearity of $D\phi(\theta_n)$,
\[
	\tau_n^{*}(\phi(\widehat\theta_n^{*})-\phi(\theta_n^{*}))-D\phi(\theta_n)\xi_n^{*}=\tau_n^{*}\big(r_n(\widehat\theta_n^{*})-r_n(\theta_n^{*})\big).
\]

Let $A_n(M)\eqdef V_n\cap \{\normin{\xi_n}\le M\}$ and $
A_n^{*}(M)\eqdef V_n\cap\left\{\norm{\xi_n^{*}}\le M\right\}\cap \left\{\tau_n^{*}\norm{\theta_n^{*}-\theta_n}\le M\right\}$, so that $\normin{\widehat\theta_n-\theta_n}\le M/\tau_n<t_n$ on $A_n(M)$, and $\normin{\widehat\theta_n^{*}-\theta_n}\le\normin{\widehat\theta_n^{*}-\theta_n^{*}}+\normin{\theta_n^{*}-\theta_n}\le2M/\tau_n^{*}<t_n$ on $A_n^{*}(M)$. Since $\ell$ is Lipschitz,
\begin{align*}
	\absin{T_n-T_n'}&\le\Lip(\ell)\varpi_n(t_n)\normin{\xi_n}\le M\Lip(\ell)\varpi_n(t_n) \qtext{on }A_n(M),
\shortintertext{and}
	\absin{T_n^{*}-T_n'^{*}}&\le\Lip(\ell)\varpi_n(t_n)\tau_n^{*}\big(\normin{\widehat\theta_n^{*}-\theta_n}+\normin{\theta_n^{*}-\theta_n}\big) \\
	&\le 3M\Lip(\ell)\varpi_n(t_n) \qtext{on }A_n^{*}(M).
\end{align*}
Hence, for any $\epsilon>0$,
\begin{align*}
	\PR{\absin{T_n-T_n'}>\epsilon\mid\C_n}&\le\PR{\normin{\xi_n}>M\mid\C_n} \\
	&\quad+\PR{M\Lip(\ell)\varpi_n(t_n)>\epsilon\mid \C_n} \qtext{a.s.},
\shortintertext{and}
	\PR{\absin{T_n^{*}-T_n'^{*}}>\epsilon\mid\C_n}&\le\PR{\normin{\xi^{*}_n}>M\mid\C_n}+\PR{\tau^{*}_n\normin{\theta^{*}_n-\theta_n}>M\mid\C_n} \\
	&\quad+\PR{3M\Lip(\ell)\varpi_n(t_n)>\epsilon\mid \C_n} \qtext{a.s.}
\end{align*}
By the assumption on the modulus $\varpi_n$, the last terms vanish a.s.\ as $n\to\infty$. Thus, by the tightness of $\seq{\normin{\xi_n}}$,
\begin{align*}
	&\limsup_{n\to\infty}\PR{\absin{T_n-T_n'}>\epsilon\mid\C_n} \\
	&\qquad\le \essinf_{M>0}\limsup_{n\to\infty}\PR{\normin{\xi_n}>M\mid\C_n}=0 \qtext{a.s.}
\end{align*}
Similarly, $\limsup_n\PR{\absin{T_n^{*}-T_n'^{*}}>\epsilon\mid\C_n}=0$ a.s. Therefore,
\[
	\absin{T_n-T_n'}\PRC{\CS}0 \qtext{and}\quad \absin{T_n^{*}-T_n'^{*}}\PRC{\CS}0 \qtext{a.s.}
\]

Finally, the result follows from Lemma \ref{lemma:nuisance} with the asymptotically equicontinuous comparison sequence $\seq{S_n'}$.
\end{proof}

Combining this result with Theorem \ref{thm:bootstrap_consistency1}, we obtain the consistency of the bootstrap for $T_n$ from conditions imposed on the linearized statistics.

\begin{corollary}
\label{corr:bootstrap_consistency3}
Suppose that the conditions of Lemma \ref{lemma:bootstrap_consistency2} hold. Then $\FD{n}{*}$ is conditionally $\KM$\-/consistent given $\CS$.
\end{corollary}

\begin{proof}
By Lemma \ref{lemma:bootstrap_consistency2}, conditions (a) and (b) of Theorem \ref{thm:bootstrap_consistency1} hold with $\seq{S_n'}$ as the comparison sequence, and $\lCZ{S_n'\mid\C_n}\to 0$ a.s.\ by the asymptotic equicontinuity of $\seq{S_n'}$.
\end{proof}

Conditions (a) and (b) of Lemma \ref{lemma:bootstrap_consistency2} involve the linear statistics $T_n'$ and $T_n'^{*}$, while merging is typically established for $\xi_n$ and $\xi_n^{*}$, for instance, by Lemma \ref{lemma:merging}. Since $T_n'$ and $T_n'^{*}$ are the images of $\xi_n$ and $\xi_n^{*}$ under the same map, such statements apply equally to the linearizations. This leads to the following counterpart of Corollary \ref{corr:bootstrap_merging}, with conditions (a) and (b) stated for $\seq{\xi_n}$ and $\seq{\xi_n^{*}}$ against a comparison sequence in $\X$.

\begin{corollary}
\label{corr:bootstrap_consistency4}
Suppose that Assumption \ref{assum:smooth_tightness} holds, $\limsup_n\normin{D\phi(\theta_n)}<\infty$ a.s., and let $\seq{Z_n}$ be a sequence of $\X$\-/valued random elements such that $\seq{S_n'}$, with $S_n'\eqdef\ell(D\phi(\theta_n)Z_n)$, is $\CS$\-/asymptotically equicontinuous. Furthermore, assume that
\begin{enumerate}[label=\tn{(\alph*)}]
	\item $\seq{\xi_n}$ merges $\CS$\-/weakly with $\seq{Z_n}$, and

	\item $\dBL{\QM_{\xi_n^{*}}^{\G_n},\QM_{Z_n}^{\C_n}}\PRC{\CS}0$ a.s.
\end{enumerate}
Then $\FD{n}{*}$ is conditionally $\KM$\-/consistent given $\CS$.
\end{corollary}

\begin{proof}
For $\omega\in\Omega$ and $x\in \X$, let $g_n(\omega,x)\eqdef\ell(D\phi(\theta_n(\omega))x)$ and $L_n\eqdef 1\vee\Lip(\ell)\normin{D\phi(\theta_n)}$. Since $\omega\mapsto D\phi(\theta_n(\omega))x$ is $\C_n$\-/measurable for every $x\in\X$, $g_n$ is a Carath\'eodory function. Moreover, $g_n$ is Lipschitz in $x$ with Lipschitz constant at most $L_n$, which is $\C_n$\-/measurable because $\normin{D\phi(\theta_n)}=\sup_{k}\normin{D\phi(\theta_n)x_k}$ for a countable dense subset $\{x_k\}$ of the unit ball of $\X$. Thus, \eqref{eq:lip_image} applies.

The conditional laws of $T_n'=g_n(\csdot,\xi_n)$, $T_n'^{*}=g_n(\csdot,\xi_n^{*})$, and $S_n'=g_n(\csdot,Z_n)$ are a.s.\ the images of the conditional laws of $\xi_n$, $\xi_n^{*}$, and $Z_n$ under $g_n$, respectively \citep[see, e.g.,][Theorem 8.5]{Kallenberg:21:FMP}. Hence, by (a) and the fact that $\limsup_n L_n<\infty$ a.s.,
\[
	\DBL{\QM_{T_n'}^{\C_n},\QM_{S_n'}^{\C_n}}\le L_n\DBL{\QM_{\xi_n}^{\C_n},\QM_{Z_n}^{\C_n}}\to 0 \qtext{a.s.}
\]
Similarly,
\[
	\DBL{\QM_{T_n'^{*}}^{\G_n},\QM_{S_n'}^{\C_n}}\le L_n\DBL{\QM_{\xi_n^{*}}^{\G_n},\QM_{Z_n}^{\C_n}} \qtext{a.s.},
\]
and, using (b), for every $\epsilon>0$, we have
\[
	\PR{\dBL{\QM_{T_n'^{*}}^{\G_n},\QM_{S_n'}^{\C_n}}>\epsilon\mid\C_n}\le\PR{\dBL{\QM_{\xi_n^{*}}^{\G_n},\QM_{Z_n}^{\C_n}}>\epsilon/L_n\mid\C_n}\to 0 \qtext{a.s.}
\]
because $\liminf_n\epsilon/L_n>0$ a.s.

Together with the asymptotic equicontinuity of $\seq{S_n'}$, the proof of Corollary \ref{corr:bootstrap_merging} provides conditions (a)--(b) of Lemma \ref{lemma:bootstrap_consistency2} for $\seq{S_n'}$. Finally, Corollary \ref{corr:bootstrap_consistency3} implies the result.
\end{proof}

\begin{remark}
\label{rem:smooth}
(1) The modulus condition in Assumption \ref{assum:smooth_tightness} holds when (i) $\rho_n^{-1}\PRC{\CS}0$ a.s.\ and (ii) $\lim_{t\downarrow 0}\limsup_n\varpi_n(t)=0$ a.s. Specifically, since $t\mapsto\varpi_n(t)$ is nondecreasing,
\begin{align}
	\label{eq:local_modulus}
	\begin{aligned}
		&\PR{\varpi_n(M/\rho_n)>\epsilon\mid\C_n} \\
		&\qquad\le\ind\{\varpi_n(1/m)>\epsilon\}+\PR{\rho_n<mM\mid\C_n} \qtext{a.s.}
	\end{aligned}
\end{align}
for every $m\in\N$. The second term on the right-hand side of \eqref{eq:local_modulus} vanishes in the limit by (i), while (ii) implies that $\limsup_n \ind\{\varpi_n(1/m)>\epsilon\}\to 0$ a.s.\ as $m\to\infty$. If $\Theta$ is convex and $\phi$ is the restriction of a continuously Fr\'echet differentiable map on an open $\Theta'\supset\Theta$, the mean value inequality \citep[see, e.g.,][(8.5.4)]{Dieudonne:69:Foundations}, applied to $z\mapsto\phi(z)-D\phi(\theta_n)z$, yields
\[
	\varpi_n(t)\le\sup\left\{\norm{D\phi(\theta)-D\phi(\theta_n)}:\theta\in\Theta,\norm{\theta-\theta_n}\le t\right\}
\]
for all $t\ge 0$ a.s. In this case, (ii) is satisfied if $\seq{\theta_n}$ is contained in a compact subset of $\Theta$ a.s., or if $\seq{\theta_n}$ is a.s.\ bounded and $D\phi$ is Lipschitz on bounded sets. Without these assumptions, (ii) holds if there are random $L<\infty$ and $t_0>0$ such that $\sup_{n\ge1}\varpi_n(t)\le Lt$ for all $t\in[0,t_0]$ a.s.

(2) No nondegeneracy condition on $D\phi(\theta_n)$ is imposed explicitly in Lemma \ref{lemma:bootstrap_consistency2}, as it is implied by the conditions of the lemma. For example, if $D\phi(\theta_n)=0$ for infinitely many $n$ on an event of positive probability, then on this event $T_n'=\ell(0)$ for the same $n$, so that condition (a) of Lemma \ref{lemma:bootstrap_consistency2} cannot hold for an asymptotically equicontinuous $\seq{S_n'}$.

(3) The proof of Lemma \ref{lemma:bootstrap_consistency2} uses only the Lipschitz property of $\ell$. Therefore, it remains valid when $\ell$ is replaced by a Carath\'eodory function $\ell_n:\Omega\times\Y\to\R$ such that $y\mapsto\ell_n(\csdot,y)$ is Lipschitz a.s.\ with a $\C_n$\-/measurable constant $L_n$, provided that the modulus condition of Assumption \ref{assum:smooth_tightness} is imposed on $L_n\varpi_n(M/\rho_n)$.
\end{remark}

The following example shows that dividing by the $\C_n$\-/measurable scale (standard deviation) of the comparison replaces the lower bound on that scale, required for the asymptotic equicontinuity of the comparison in Lemma \ref{lemma:bootstrap_consistency2} and Corollary \ref{corr:bootstrap_consistency4}, with a stronger modulus condition.

\begin{example}
\label{ex:rescaling}
Take $\X=\R^{d}$, $\Y=\R$, and $\ell=\mathrm{id}$. Let $Z_n\mid \C_n\sim\nD{0}{\Sigma_n}$, and let $g_n\eqdef D\phi(\theta_n)$ and $S_n'\eqdef g_n^{\top}Z_n$, so that $\QM_{S_n'}^{\C_n}=\nD{0}{\sigma_n^{2}}$ with $\sigma_n^{2}\eqdef g_n^{\top}\Sigma_ng_n>0$ a.s. The sequence $\seq{S_n'}$ is $\CS$\-/asymptotically equicontinuous if $\liminf_n\sigma_n>0$ a.s.\ because $\lC{\epsilon,S_n'\mid\C_n}\le\epsilon/(\sigma_n\sqrt{2\pi})$ a.s. Alternatively, by Lemma \ref{lemma:scale_inv}, conditions (a) and (b) and the conclusion of the lemma are unaffected by dividing $T_n'$, $T_n'^{*}$ and $S_n'$ by the $\C_n$\-/measurable $\sigma_n$, and, by Remark \ref{rem:smooth}(3) with $\ell_n=\ell/\sigma_n$, the lemma applies to $T_n/\sigma_n$ and $T_n^{*}/\sigma_n$ with the modulus condition imposed on $\varpi_n(M/\rho_n)/\sigma_n$.
\end{example}

\subsection{Degenerate first derivative}
\label{sec:second_order}

When the first derivative of $\phi$ vanishes at $\theta_n$, the linearized statistics in \eqref{eq:linear_stat} are degenerate, the expansion requires a second-order extension, and the bootstrap statistic must be corrected by subtracting the estimated linear term. We assume that $\Theta$ is a convex Borel subset of $\X$ and $\phi$ is the restriction of a twice continuously Fr\'echet differentiable map on an open $\Theta'\supset\Theta$. For $\theta\in\Theta'$, the derivative of $\phi$ at $\theta$ is a bounded linear map $D\phi(\theta):\X\to\Y$, assumed to vanish at $\theta=\theta_n$ a.s.\ for each $n\ge 1$, and its second derivative is identified with a bounded symmetric bilinear map $D^{2}\phi(\theta):\X\times\X\to\Y$. The remaining notation is that of Section \ref{sec:linearization}. Consider the modified statistics
\begin{align*}
	\widecheck{T}_n&\eqdef\ell\big(\tau_n^{2}(\phi(\widehat\theta_n)-\phi(\theta_n))\big) \qtext{and} \\
	\widecheck{T}_n^{*}&\eqdef\ell\big(\tau_n^{*2}(\phi(\widehat{\theta}_n^{*})-\phi(\theta_n^{*})-D\phi(\theta_n^{*})(\widehat{\theta}_n^{*}-\theta_n^{*}))\big).
\end{align*}
Their quadratic approximations are given by
\[
	\widecheck{T}_n'\eqdef\ell\left(\frac{1}{2}D^{2}\phi(\theta_n)[\xi_n,\xi_n]\right) \qtext{and}\quad	\widecheck{T}_n'^{*}\eqdef\ell\left(\frac{1}{2}D^{2}\phi(\theta_n)[\xi_n^{*},\xi_n^{*}]\right),
\]
respectively. For $t\ge 0$, let
\[
	\widecheck{\varpi}_n(t)\eqdef\sup\left\{\norm{D^{2}\phi(\theta)-D^{2}\phi(\theta_n)}:\theta\in\Theta, \norm{\theta-\theta_n}\le t\right\}
\]
denote the local modulus of continuity of $D^{2}\phi$ at $\theta_n$.\footnote{
	Since $\Theta$ is convex and $D^{2}\phi$ is continuous, the supremum is unchanged if the weak inequality is replaced by the strict one (for $t>0$), and the resulting function of $(\theta_n,t)$ is lower semicontinuous. Hence, $\widecheck{\varpi}_n(R)$ is $\F$\-/measurable for every $\sigma$-algebra $\F$ with respect to which $\theta_n$ and $R$ are measurable.
}

\begin{corollary}
\label{corr:bootstrap_consistency5}
Suppose that Assumption \ref{assum:smooth_tightness} holds with $\widecheck{\varpi}_n$ in place of $\varpi_n$, and $\seq{\widecheck{S}_n'}$ is a $\CS$\-/asymptotically equicontinuous sequence of random variables. Furthermore, assume that
\begin{enumerate}[label=\tn{(\alph*)}]
	\item $\dK{\widecheck{T}_n',\widecheck{S}_n'\mid\C_n}\to0$ a.s., and

	\item $\dK{\widecheck{T}_n'^{*},\widecheck{S}_n'\mid\G_n,\C_n}\PRC{\CS}0$ a.s.
\end{enumerate}
Then the conclusion of Theorem \ref{thm:bootstrap_consistency1} holds for $\widecheck T_n$ and $\widecheck T_n^{*}$.
\end{corollary}

\begin{proof}
Let $t_n$, $A_n(M)$ and $A_n^{*}(M)$ be as in the proof of Lemma \ref{lemma:bootstrap_consistency2}. Define
\[
	\psi_n(z)\eqdef\phi(z)-\frac{1}{2}D^{2}\phi(\theta_n)[z-\theta_n,z-\theta_n],
\]
so that $D^{2}\psi_n(z)=D^{2}\phi(z)-D^{2}\phi(\theta_n)$ by the symmetry of $D^{2}\phi(\theta_n)$ \citep[see, e.g.,][(8.12.2)]{Dieudonne:69:Foundations}. For $x,y\in\Theta$ and $h\eqdef y-x$, we have
\begin{align}
	\label{eq:psi_eq}
	\begin{aligned}
		&\psi_n(y)-\psi_n(x)-D\psi_n(x)h \\
		&\qquad=\phi(y)-\phi(x)-D\phi(x)h-\frac{1}{2}D^{2}\phi(\theta_n)[h,h].
	\end{aligned}
\end{align}
Applying the mean value inequality to $z\mapsto\psi_n(z)-D\psi_n(x)z$ along the segment joining $x$ and $y$, and then to $D\psi_n$ along the segment joining $x$ and $x+sh$ for each $s\in[0,1]$, we bound the norm of the right-hand side of \eqref{eq:psi_eq} by
\begin{equation}
	\label{eq:psi_bound}
	\norm{h}^{2}\sup_{0\le s\le 1}\norm{D^{2}\phi(x+sh)-D^{2}\phi(\theta_n)}.
\end{equation}

We apply the bound in \eqref{eq:psi_bound} with $(x,y)=(\theta_n,\widehat\theta_n)$ and $(x,y)=(\theta_n^{*},\widehat\theta_n^{*})$. Since $\Theta$ is convex, both segments lie in the ball of radius $t_n$ around $\theta_n$ on $A_n(M)$ and $A_n^{*}(M)$, respectively, and $D\phi(\theta_n)=0$ a.s., so that
\begin{align*}
	\absin{\widecheck{T}_n-\widecheck {T}_n'}&\le M^{2}\Lip(\ell)\widecheck{\varpi}_n(t_n) \qtext{on }A_n(M),
\shortintertext{and}
	\absin{\widecheck{T}_n^{*}-\widecheck{T}_n'^{*}}&\le M^{2}\Lip(\ell)\widecheck{\varpi}_n(t_n) \qtext{on }A_n^{*}(M).
\end{align*}
Following the proof of Lemma \ref{lemma:bootstrap_consistency2}, we get $\dK{\widecheck T_n,\widecheck S_n'\mid\C_n}\to0$ and $\dK{\widecheck T_n^{*},\widecheck S_n'\mid\G_n,\C_n}\PRC{\CS}0$ a.s., and the result follows from Theorem \ref{thm:bootstrap_consistency1}.
\end{proof}

The sufficient conditions (i) and (ii) of Remark \ref{rem:smooth}(1) apply to $\widecheck{\varpi}_n$ as well. For instance, (ii) is satisfied if $\seq{\theta_n}$ is contained in a compact subset of $\Theta$ a.s., or if $\seq{\theta_n}$ is a.s.\ bounded and $D^{2}\phi$ is Lipschitz on bounded sets. Similarly, Remark \ref{rem:smooth}(3) extends to Corollary \ref{corr:bootstrap_consistency5}, with the modulus condition imposed on $L_n\widecheck{\varpi}_n(M/\rho_n)$. No nondegeneracy condition on $D^{2}\phi(\theta_n)$ is imposed explicitly as it is implied by condition (a) of the corollary.

\subsection{Smoothness in the nuisance}
\label{sec:nuisance_smooth}

The linearization of Section \ref{sec:linearization} handles the dependence of the statistic on the estimator. We now consider the case of the estimated nuisance discussed in Section \ref{sec:setup}. Write
\[
	T_n=\Psi_n(\zeta_n,\vartheta_n) \qtext{and}\quad \widetilde{T}_n=\Psi_n(\zeta_n,\widehat\vartheta_n),
\]
where $\zeta_n$ is a $\Z$\-/valued data-based random element, $\vartheta_n\in\Theta_{\vartheta}$ a.s.\ for a convex Borel set $\Theta_{\vartheta}\subset\Y$, $\widehat\vartheta_n\in\Theta_{\vartheta}$ a.s.\ is an estimator of $\vartheta_n$, and for an open $\Theta_{\vartheta}'\supset \Theta_{\vartheta}$, $\Psi_n:\Z\times\Theta_{\vartheta}'\to\R$ is measurable. The bootstrap counterparts are given by
\[
	T_n^{*}=\Psi_n(\zeta_n^{*},\vartheta_n) \qtext{and}\quad \widetilde{T}_n^{*}=\Psi_n(\zeta_n^{*},\widehat\vartheta_n^{*}),
\]
where $\zeta_n^{*}$ and $\widehat{\vartheta}_n^{*}$ are the resampling versions of $\zeta_n$ and $\widehat{\vartheta}_n$, respectively, with $\widehat{\vartheta}_n^{*}\in\Theta_{\vartheta}$ a.s.

By Lemma \ref{lemma:nuisance}, the feasible pair satisfies conditions (a) and (b) of Theorem \ref{thm:bootstrap_consistency1}, provided that the infeasible pair satisfies these conditions with a $\CS$\-/asymptotically equicontinuous comparison sequence $\seq{S_n}$, and
\[
	\absin{\widetilde{T}_n-T_n}\PRC{\CS}0 \qtext{and}\quad \absin{\widetilde{T}_n^{*}-T_n^{*}}\PRC{\CS}0 \qtext{a.s.}
\]
For example, when $\Y=\Z=\R$ and $\Psi_n(z,\vartheta)=z-\vartheta$, the latter follows from the consistency of $\widehat\vartheta_n$ and $\widehat\vartheta_n^{*}$. For $\Theta_{\vartheta}=(0,\infty)$ and $\Psi_n(z,\vartheta)=z/\vartheta$, we have $\widetilde{T}_n=a_nT_n$ with $a_n\eqdef\vartheta_n/\widehat\vartheta_n$. Suppose that $\PR{S_n=0\mid\C_n}=0$ a.s. The argument of Lemma \ref{lemma:nuisance}, with intervals $[x,(1+2\epsilon)x]$, $x>0$, and $[(1+2\epsilon)x,x]$, $x<0$, in place of $[x,x+\epsilon]$, implies that for $\epsilon\in(0,1/2]$,
\begin{align*}
	\DK{a_nT_n,S_n\mid\C_n}&\le\DK{T_n,S_n\mid\C_n}+\LC{2\epsilon,\ln\absin{S_n}\mid\C_n} \\
	&\quad+\PR{\abs{a_n-1}>\epsilon\mid\C_n} \qtext{a.s.},
\end{align*}
because $\ln\absin{\csdot}$ maps these intervals onto intervals of length $\ln(1+2\epsilon)\le2\epsilon$. The bootstrap pair is handled similarly, with $a_n^{*}\eqdef\vartheta_n/\widehat\vartheta_n^{*}$ and conditioning on $\G_n$ in the first and last terms. Since $\ln\absin{\sigma G}=\ln\sigma+\ln\absin{G}$ for $G\sim\nD{0}{1}$, the L\'evy concentration function of $\ln\absin{S_n}$ is invariant under rescaling of $S_n$, so that this bound applies to Gaussian comparisons with degenerating scales. The following lemma covers a more general case.

\begin{lemma}
\label{lemma:nuisance_cont}
Suppose that $\vartheta_n\in C$ a.s.\ for each $n\ge 1$ and a compact $C\subset\Theta_{\vartheta}$, that $\normin{\widehat\vartheta_n-\vartheta_n}\PRC{\CS}0$ and $\normin{\widehat\vartheta_n^{*}-\vartheta_n}\PRC{\CS}0$ a.s., and that $\seq{\zeta_n}$ and $\seq{\zeta_n^{*}}$ are $\CS$\-/asymptotically tight. If for every $L\in\K_{\Z}$, the family $\{\Psi_n(z,\csdot):n\ge 1,\ z\in L\}$ is equicontinuous at each point of $C$, then $\absin{\widetilde{T}_n-T_n}\PRC{\CS}0$ and $\absin{\widetilde{T}_n^{*}-T_n^{*}}\PRC{\CS}0$ a.s.
\end{lemma}

\begin{proof}
Fix $\epsilon>0$ and $L\in\K_{\Z}$. By equicontinuity, for each $x\in C$ there is $\rho_x>0$ such that $\absin{\Psi_n(z,x)-\Psi_n(z,y)}<\epsilon$ for all $n\ge1$, $z\in L$, and $y\in\Theta_{\vartheta}'$ with $\normin{x-y}<\rho_x$. Cover the compact set $C$ by finitely many balls $B(x_i,\rho_{x_i}/2)$, $1\le i\le N$, and let $\delta_L\eqdef\min_{i\le N}\rho_{x_i}/2$. If $x\in C$ and $\normin{x-y}<\delta_L$, choose $i$ with $\normin{x-x_i}<\rho_{x_i}/2$; then $\normin{y-x_i}<\rho_{x_i}$, so that $\absin{\Psi_n(z,x)-\Psi_n(z,y)}<2\epsilon$ for all $n\ge1$ and $z\in L$. Since $\vartheta_n\in C$ a.s.,
\begin{align*}
	\PR{\absin{\widetilde{T}_n-T_n}\ge 2\epsilon\mid\C_n}&\le\PR{\normin{\widehat\vartheta_n-\vartheta_n}\ge\delta_L\mid\C_n} \\
	&\quad+\PR{\zeta_n\notin L\mid\C_n} \qtext{a.s.}
\end{align*}
The first term vanishes a.s.\ as $n\to\infty$, and taking the essential infimum over $L\in\K_{\Z}$ yields $\absin{\widetilde{T}_n-T_n}\PRC{\CS}0$ a.s. The same argument with $\zeta_n^{*}$ and $\widehat\vartheta_n^{*}$ in place of $\zeta_n$ and $\widehat\vartheta_n$ applies to the bootstrap pair.
\end{proof}

Finally, when $\Psi_n$ is differentiable in $\vartheta$, the requirement reduces to conditions on the estimator and the derivative of $\Psi_n$. Assume that for every $z\in\Z$, $\vartheta\mapsto\Psi_n(z,\vartheta)$ is continuously Fr\'echet differentiable on $\Theta_{\vartheta}'$, and its derivative $\partial_{\vartheta}\Psi_n(z,\vartheta)$ is a bounded linear functional on $\Y$. Define $V_n'\eqdef\{(\vartheta_n,\widehat\vartheta_n,\widehat{\vartheta}_n^{*})\in\Theta_{\vartheta}^{3}\}$,
\begin{align*}
	\widebar{G}_n&\eqdef\sup_{t\in[0,1]\cap\Q}\norm{\partial_{\vartheta}\Psi_n\big(\zeta_n,\vartheta_n+t(\widehat{\vartheta}_n-\vartheta_n)\big)} \qtext{and} \\
	\widebar{G}_n^{*}&\eqdef\sup_{t\in[0,1]\cap\Q}\norm{\partial_{\vartheta}\Psi_n\big(\zeta_n^{*},\vartheta_n+t(\widehat{\vartheta}_n^{*}-\vartheta_n)\big)}.
\end{align*}

\begin{lemma}
\label{lemma:nuisance_smooth}
Suppose that $\normin{\widehat{\vartheta}_n-\vartheta_n}\PRC{\CS}0$ and $\normin{\widehat{\vartheta}_n^{*}-\vartheta_n}\PRC{\CS}0$ a.s., and $\seq{\widebar{G}_n}$ and $\seq{\widebar{G}_n^{*}}$ are $\CS$\-/asymptotically tight. Then $\absin{\widetilde{T}_n-T_n}\PRC{\CS}0$ and $\absin{\widetilde{T}_n^{*}-T_n^{*}}\PRC{\CS}0$ a.s. Consequently, if the infeasible pair $(T_n,T_n^{*})$ satisfies conditions \tn{(a)} and \tn{(b)} of Theorem \ref{thm:bootstrap_consistency1} with a $\CS$\-/asymptotically equicontinuous comparison sequence $\seq{S_n}$, the feasible bootstrap is conditionally $\KM$\-/consistent given $\CS$.
\end{lemma}

\begin{proof}
Fix $M>0$ and let $A_n(M)\eqdef V_n'\cap\{\widebar{G}_n\le M\}$ and $A_n^{*}(M)\eqdef V_n'\cap\{\widebar{G}_n^{*}\le M\}$. On $V_n'$, the segments joining $\vartheta_n$ to $\widehat\vartheta_n$ and to $\widehat\vartheta_n^{*}$ lie in $\Theta_{\vartheta}$ by the convexity of $\Theta_{\vartheta}$. Hence, by the mean value inequality applied to $\vartheta\mapsto\Psi_n(\zeta_n,\vartheta)$ and $\vartheta\mapsto\Psi_n(\zeta_n^{*},\vartheta)$ along these segments,
\begin{align*}
	\absin{\widetilde{T}_n-T_n}&\le M\normin{\widehat{\vartheta}_n-\vartheta_n} \qtext{on $A_n(M)$, and} \\ \absin{\widetilde{T}_n^{*}-T_n^{*}}&\le M\normin{\widehat\vartheta_n^{*}-\vartheta_n} \qtext{on }A_n^{*}(M).
\end{align*}
Thus, for $\delta>0$,
\begin{align}
	\label{eq:nuisance_split}
	\begin{aligned}
		\PR{\absin{\widetilde{T}_n-T_n}>\delta\mid\C_n}&\le\PR{\normin{\widehat{\vartheta}_n-\vartheta_n}>\delta/M\mid\C_n} \\
		&\quad+\PR{\widebar{G}_n>M\mid\C_n} \qtext{a.s.}
	\end{aligned}
\end{align}
The first term on the right-hand side of \eqref{eq:nuisance_split} vanishes a.s.\ as $n\to\infty$, and $\seq{\widebar{G}_n}$ is $\CS$\-/asymptotically tight. Therefore, $\absin{\widetilde{T}_n-T_n}\PRC{\CS}0$ a.s. The bootstrap pair is handled similarly, with $\widebar{G}_n$ and $A_n$ replaced by $\widebar G_n^{*}$ and $A_n^{*}$. Consequently, the result follows from Theorem \ref{thm:bootstrap_consistency1}, where $\seq{S_n}$ is the comparison sequence, and conditions (a)--(b) hold by Lemma \ref{lemma:nuisance}.
\end{proof}

\begin{example}
Take $\Z=\R^{d}$ and let $\Y$ be the space of symmetric $d\times d$ matrices, in which the cone $\Theta_{\vartheta}$ of positive definite matrices is open. The Mahalanobis norm and its square,
\[
	\Psi_n(z,\vartheta)=\normin{\vartheta^{-1/2}z} \qtext{and}\quad \Psi_n(z,\vartheta)=z^{\top}\vartheta^{-1}z,
\]
are continuously differentiable in $\vartheta$, with $\normin{\partial_{\vartheta}\Psi_n(z,\vartheta)}\le \normin{\vartheta^{-1}}^{3/2}\normin{z}$ in the first case and $\normin{\partial_{\vartheta}\Psi_n(z,\vartheta)}\le\normin{\vartheta^{-1}}^{2}\normin{z}^{2}$ in the second. Since $\normin{\vartheta^{-1}}=1/\lambda_{\min}(\vartheta)$ and, by Weyl's inequality, $\lambda_{\min}(\vartheta)\ge\lambda_{\min}(\vartheta_n)-\normin{\vartheta-\vartheta_n}$, we have $\normin{\vartheta^{-1}}\le2\normin{\vartheta_n^{-1}}$ on the segment in the definition of $\widebar G_n$ whenever $\normin{\widehat\vartheta_n-\vartheta_n}\le\lambda_{\min}(\vartheta_n)/2$. Hence, for the squared norm,
\[
	\widebar G_n\le4\normin{\vartheta_n^{-1}}^{2}\normin{\zeta_n}^{2} \qtext{on } \big\{\normin{\widehat\vartheta_n-\vartheta_n}\le\lambda_{\min}(\vartheta_n)/2\big\},
\]
and similarly for the norm, with $2^{3/2}\normin{\vartheta_n^{-1}}^{3/2}\normin{\zeta_n}$, and for $\widebar G_n^{*}$, with $\zeta_n^{*}$ and $\widehat\vartheta_n^{*}$. Since $\normin{\widehat\vartheta_n-\vartheta_n}\PRC{\CS}0$ and $\normin{\widehat\vartheta_n^{*}-\vartheta_n}\PRC{\CS}0$ a.s., $\seq{\widebar{G}_n}$ and $\seq{\widebar{G}_n^{*}}$ are $\CS$\-/asymptotically tight if $\seq{\zeta_n}$, $\seq{\zeta_n^{*}}$ and $\seq{\normin{\vartheta_n^{-1}}}$ are $\CS$\-/asymptotically tight.

For instance, if $\vartheta_n^{-1/2}\zeta_n\to\nD{0}{I_d}$, $\CS$\-/weakly w.r.t.\ $\H_{\varnothing}$, then the infeasible quadratic form $\zeta_n^{\top}\vartheta_n^{-1}\zeta_n$ converges $\CS$\-/weakly to the fixed law $\chi^{2}_{d}$. In this case, Corollary \ref{corr:bootstrap_consistency1} applies to the infeasible pair with $S\sim\chi^{2}_{d}$ independent of everything else, and Lemma \ref{lemma:nuisance_smooth} transfers the conclusion to the feasible bootstrap.
\end{example}

\begin{remark}
\label{rem:sfm_varying}
The proofs in this section bound the relevant differences separately for each $n$, using only the norms of $\X$ and $\Y$, the linearity of the derivative or the mean value inequality, and the Lipschitz property of $\ell$. Consequently, all the results except Lemma \ref{lemma:nuisance_cont} remain valid when $\X$, $\Y$, $\Z$, the parameter sets, and the maps $\phi$ and $\ell$ depend on $n$, with $\Lip(\ell)$ replaced by $\Lip(\ell_n)$ and $\limsup_n\Lip(\ell_n)<\infty$, a condition relaxed by Remark \ref{rem:smooth}(3). This covers high-dimensional settings, in which the growth of the dimension is absorbed by the normalization $\tau_n$ and by the Lipschitz constant of $\ell_n$ (see the example of Section \ref{sec:example_corr}). When $\Theta$ and $\phi$ depend on $n$, the sufficient conditions for (ii) in Remark \ref{rem:smooth}(1) and in Section \ref{sec:second_order} that rely on compactness or boundedness of $\seq{\theta_n}$ are not available, and the modulus condition is verified directly, as in Section \ref{sec:example_corr}, or through the last condition of Remark \ref{rem:smooth}(1).
\end{remark}

\section{Examples}
\label{sec:examples}

\subsection{Maxima in fixed dimension}
\label{sec:example_fd}

This example is inspired by the wobbling maxima case of \citet[Section 5]{DAristotile/Diaconis/Freedman:88} and chosen because the conditional laws of the statistic do not converge weakly, so that the comparison must vary with the sample size. Moreover, no $\C_n$\-/measurable affine transformation leads to convergence to a nondegenerate fixed kernel.

Let $\seq{(m_n,\rho_n)}$ be i.i.d.\ $\R^{2}\times \R$\-/valued random elements such that $\rho_n$ is nondegenerate and $\supp(\rho_n)\subset[-\widebar{\rho},\widebar{\rho}]$ for some $\widebar{\rho}\in (0,1)$. Take $\C_n=\sigma(m_n,\rho_n)$ and let $R_n$ denote the correlation matrix with off-diagonal entry $\rho_n$. For each $n\ge 1$, let $X_{n,1},\ldots,X_{n,n}$ be a martingale difference sequence in $\R^{2}$ with respect to a filtration $\F_{n,1}\subset \cdots \subset \F_{n,n}$ such that $\F_{n,0}\eqdef\C_n\subset\F_{n,1}$, and for $1\le i\le n$, $\E[X_{n,i}X_{n,i}^{\top}\mid\F_{n,i-1}]=I_2$ and $\E[\normin{X_{n,i}}^3\mid\C_n]\le K<\infty$ a.s. Let $R_n^{1/2}$ be the symmetric square root of $R_n$. The observed sample consists of $Y_{n,i}$, $i=1,\ldots,n$, where
\[
	Y_{n,i}\eqdef m_n+\zeta_{n,i}, \quad \zeta_{n,i}\eqdef R_n^{1/2}X_{n,i}.
\]
In this case, $\E[Y_{n,i}\mid\C_n]=m_n$ and $\Var(Y_{n,i}\mid\C_n)=R_n$ a.s., that is, the shock affects the location of the coordinates and their correlation. The parameter of interest is $\theta_n=m_n$.

Let $Y_n^{*}$ be drawn i.i.d.\ from the empirical distribution of $Y_n$ and consider the estimators $\widehat{\theta}_n=\theta_n^{*}=\widebar{Y}_n$. Define $\xi_n\eqdef\sqrt{n}(\widebar{Y}_n-\theta_n)$, $\xi_n^{*}\eqdef\sqrt{n}(\widebar{Y}_n^{*}-\widebar{Y}_n)$, and $Z_n\eqdef R_n^{1/2}Z$, where $Z$ is a standard normal random vector in $\R^{2}$, independent of everything else. We consider the statistics
\[
	T_n=\max\{\xi_{n,1},\xi_{n,2}\} \qtext{and}\quad T_n^{*}=\max\{\xi_{n,1}^{*},\xi_{n,2}^{*}\},
\]
and the comparison $S_n=\max\{Z_{n,1},Z_{n,2}\}$.

We check the conditions of Theorem \ref{thm:bootstrap_consistency1}. Note that $\lCZ{S_n\mid\C_n}=0$ a.s. For condition (a), the vectors $\zeta_{n,i}$, $1\le i\le n$, form a martingale difference sequence with respect to the same filtration. Additionally, $\E[\zeta_{n,i}\zeta_{n,i}^{\top}\mid \F_{n,i-1}]=R_n$ a.s., which is $\C_n$\-/measurable with the smallest eigenvalue $1-\absin{\rho_n}\ge1-\widebar{\rho}$ a.s., and $\E[\norm{\zeta_{n,i}}_{\infty}^3\mid\C_n]\le K$ a.s. Since $\sum_{i=1}^n \zeta_{n,i}=\sqrt n\xi_n$ and $\sqrt nZ_n$ is conditionally normal given $\C_n$ with the same conditional covariance matrix $nR_n$, the Berry--Esseen bound of \citet[Theorem 3.1 and the remark following it]{Kojevnikov/Song:22} yields $\Delta_n\le C_1n^{-1/4}$ a.s., where $C_1>0$ depends only on $K$ and $\widebar{\rho}$.

For condition (b), note that $R_n$ is a.s.\ invertible, so that $X_{n,i}=R_n^{-1/2}(Y_{n,i}-m_n)$ a.s., and conditional laws given $\G_n$ coincide a.s.\ with those given $\C_n\vee\sigma(X_{n,1},\ldots,X_{n,n})$. Let $X_n^{*}$ denote the bootstrap analog of $X_{n,1},\ldots,X_{n,n}$ corresponding to $Y_n^{*}$. Define $V_n^{*}\eqdef\sqrt n(\widebar{X}_n^{*}-\widebar{X}_n)$, so that $\xi_n^{*}=R_n^{1/2}V_n^{*}$, and let $\QM_n^{*}$ be the conditional law of $V_n^{*}$ given $\G_n$. Conditionally on $\G_n$, the elements of $X_n^{*}$ are i.i.d.\ draws from the empirical distribution of $X_{n,1},\ldots,X_{n,n}$. Therefore, by the identity in Remark \ref{rem:merging1}, $\QM_n^{*}$ is a.s.\ the law of a normalized sum of $n$ i.i.d.\ random vectors with mean zero and covariance matrix $\widehat\Sigma_n\eqdef n^{-1}\sum_{i=1}^n X_{n,i}X_{n,i}^{\top}-\widebar{X}_n\widebar{X}_n^{\top}$.

Let $\GM_{\Sigma}$ denote the probability kernel corresponding to $\nD{0}{\Sigma}$. Since $\{T_n^{*}\le t\}$ and $\{S_n\le t\}$ are the events that $V_n^{*}$ and $Z$ lie in the intersection of two half-planes determined by $R_n$ and $t$, by the identity in Remark \ref{rem:merging1},
\[
	\Delta_n^{*}\le\sup_{A\in \A}\absin{\QM_n^{*}(\csdot,A)-\GM_{\widehat\Sigma_n}(\csdot,A)}+\DTV{\GM_{\widehat{\Sigma}_n},\GM_{I_2}} \qtext{a.s.},
\]
where $\A\eqdef\{A\in \B(\R^2): A\text{ is convex}\}$ and $\dTV{\csdot,\csdot}$ denotes the total variation distance. Let $\Xi_n\eqdef\{\normin{\widehat{\Sigma}_n-I_2}_{\mathrm{F}}\le1/2\}$, where $\normin{\csdot}_{\mathrm{F}}$ is the Frobenius norm, so that the smallest eigenvalue of $\widehat{\Sigma}_n$ is at least $1/2$ on $\Xi_n$. Since linear bijections preserve convexity, the Berry--Esseen bound for convex sets \citep[see, e.g.,][Corollary~17.2]{Bhattacharya/Rao:76:NormalApprox}, applied $\omega$-wise, yields
\[
	\sup_{A\in \A}\absin{\QM_n^{*}(\csdot,A)-\GM_{\widehat\Sigma_n}(\csdot,A)}\le C_2M_n n^{-1/2} \qtext{a.s. on } \Xi_n,
\]
where $M_n\eqdef n^{-1}\sum_{i=1}^n \normin{X_{n,i}}^{3}$, and $C_2>0$ is a universal constant. Moreover, by Pinsker's inequality,
\[
	\dTV{\GM_{\widehat{\Sigma}_n},\GM_{I_2}}\le\normin{\widehat\Sigma_n-I_2}_{\mathrm{F}}/2 \qtext{a.s.\ on } \Xi_n.
\]
Hence, for every $\epsilon>0$,
\begin{align}
	\label{eq:ex1}
	\begin{aligned}
		\PR{\Delta_n^{*}>\epsilon\mid\C_n}&\le\PR{\normin{\widehat{\Sigma}_n-I_2}_{\mathrm{F}}>\epsilon\wedge (1/2)\mid\C_n} \\
		&\quad+\PR{C_2 M_n/\sqrt n>\epsilon/2\mid\C_n} \qtext{a.s.}
	\end{aligned}
\end{align}

Let $A_n\eqdef n^{-1}\sum_{i=1}^n(X_{n,i}X_{n,i}^{\top}-I_2)$, so that $\widehat\Sigma_n-I_2=A_n-\widebar X_n\widebar X_n^{\top}$. The entries of $X_{n,i}X_{n,i}^{\top}-I_2$, $1\le i\le n$, are martingale differences with respect to the same filtration, and their conditional absolute moments of order $3/2$ given $\C_n$ are bounded by $2K$ a.s. The argument of Lemma 1 of \cite{Chatterji:69} yields $\E[\absin{[A_n]_{kl}}^{3/2}\mid\C_n]\le 2\sqrt2Kn^{-1/2}$ a.s.\ for $k,l=1,2$. Since $\normin{A_n}_{\mathrm{F}}^{3/2}\le\sum_{k,l}\absin{[A_n]_{kl}}^{3/2}$, we get
\[
	\E[\normin{A_n}_{\mathrm{F}}\mid\C_n]\le\big(\E[\normin{A_n}_{\mathrm{F}}^{3/2}\mid\C_n]\big)^{2/3}\le 2^{7/3}K^{2/3}n^{-1/3} \qtext{a.s.}
\]
by the conditional Jensen inequality. Moreover, $\E[\normin{\widebar X_n\widebar X_n^{\top}}_{\mathrm{F}}\mid\C_n]=\E[\normin{\widebar X_n}^2\mid\C_n]=2/n$ a.s., and $\E[M_n\mid\C_n]\le K$ a.s. Thus, the right-hand side of \eqref{eq:ex1} vanishes a.s.\ as $n\to\infty$ by the conditional Markov inequality.

Consequently, the bootstrap is conditionally $\KM$\-/consistent given $\CS$ by Theorem \ref{thm:bootstrap_consistency1}. Finally, the simultaneous lower confidence bounds
\[
	\theta_{n,j}\ge\widebar{Y}_{n,j}- n^{-1/2}c_n^{*}(1-\alpha), \quad j=1,2,
\]
have asymptotically correct conditional coverage, uniformly over $\alpha\in (0,1)$.

\begin{remark*}
(1) Since each coordinate of $Z_n$ is standard normal given $\C_n$, $\PR{\absin{S_n}>M\mid\C_n}\le2\PR{\absin{G}>M}$ a.s.\ for each $M>0$ and $G\sim\nD{0}{1}$, so that $\seq{S_n}$ is $\CS$\-/asymptotically tight. Hence, by condition (a) and Lemma \ref{lemma:tight_conv}, $\seq{T_n}$ merges $\CS$\-/weakly with $\seq{S_n}$. Moreover, $T_n$ converges unconditionally. Since the Kolmogorov distance between the laws of $T_n$ and $S_n$ is bounded by $\E[\Delta_n]$ and $\seq{S_n}$ are identically distributed, condition (a) yields $T_n\Rightarrow S_1$.

Corollary \ref{corr:bootstrap_merging} provides an alternative route to consistency that does not rely on Berry--Esseen bounds. Suppose that $X_{n,i}=X_i$, where $\seq{X_i}$ are i.i.d.\ with mean zero and identity covariance matrix, independent of $\seq{(m_n,\rho_n)}$. In this case, by the identity in Remark \ref{rem:merging1} and \eqref{eq:lip_image}, applied to the $1$\-/Lipschitz map $x\mapsto\max_{1\le j\le2}[R_n^{1/2}x]_j$, conditions (a) and (b) of Corollary \ref{corr:bootstrap_merging} follow from the central limit theorem for $\sqrt n\widebar X_n$ and its bootstrap analog \citep{Bickel/Freedman:81} together with Lemma \ref{lemma:as_to_cond}(ii), while $\LC{\epsilon,S_n\mid\C_n}\le2\epsilon/\sqrt{2\pi}$ a.s. Note that this route requires only finite second moments.

(2) The conditional laws of $T_n$ do not converge to a nondegenerate kernel under $\C_n$\-/measurable affine transformations. For $\rho\in[-\widebar{\rho},\widebar{\rho}]$ and independent $G,H\sim\nD{0}{1}$, let $\L_{\rho}$ denote the law of $\sqrt{(1+\rho)/2}\,G+\sqrt{(1-\rho)/2}\,\absin{H}$. Since $2\cdot\max\{z_1,z_2\}=(z_1+z_2)+\absin{z_1-z_2}$, and $Z_{n,1}+Z_{n,2}$ and $Z_{n,1}-Z_{n,2}$ are conditionally uncorrelated and jointly normal given $\C_n$ with conditional variances $2(1\pm\rho_n)$, we have $\QM_{S_n}^{\C_n}=\L_{\rho_n}$ a.s. Each $\L_{\rho}$ is nondegenerate, $\rho\mapsto\L_{\rho}$ is weakly continuous, and the skewness of $\L_{\rho}$ is given by
\[
	\gamma(\rho)\eqdef\frac{4-\pi}{2}\left(\frac{1-\rho}{\pi-(1-\rho)}\right)^{3/2},
\]
which is strictly decreasing on $[-\widebar{\rho},\widebar{\rho}]$.

Let $W_n\eqdef\alpha_nT_n+\beta_n$ for $\C_n$\-/measurable $\alpha_n\ne0$ and $\beta_n$, and suppose that $\QM_{W_n}^{\C_n}\Rightarrow\QM$. Since the $\seq{\C_n}$ are independent, $\QM$ equals a.s.\ a fixed law $\L$ by Lemma \ref{lemma:tail_join}(c), and we assume that $\L$ is nondegenerate. Fix $\omega\in\Omega$ for which $\QM_{W_n}^{\C_n}(\omega,\csdot)\Rightarrow\L$ and the merging in (1) holds, and an infinite $N'\subset\N$ such that $\rho_n(\omega)\to\rho'\in[-\widebar{\rho},\widebar{\rho}]$ and the sign of $\alpha_n(\omega)$ is constant along $n\in N'$. Then $\QM_{T_n}^{\C_n}(\omega,\csdot)\Rightarrow\L_{\rho'}$ along $n\in N'$ because $\QM_{T_n}^{\C_n}\mrg\QM_{S_n}^{\C_n}=\L_{\rho_n}$, and the convergence of types theorem \citep[Theorem 14.2]{Billingsley:95:Probability}, applied to $\pm W_n$, implies that $\L$ is the image of $\L_{\rho'}$ under an affine map with nonzero slope. In particular, $\L$ has a finite third moment, and the skewness $\gamma'$ of $\L$ satisfies $\absin{\gamma'}=\gamma(\rho')$. Hence, every limit point of $\seq{\rho_n(\omega)}$ equals $\gamma^{-1}(\absin{\gamma'})$, that is, $\rho_n\to\gamma^{-1}(\absin{\gamma'})$ a.s., which is impossible because $\seq{\rho_n}$ are identically distributed and nondegenerate. For $\alpha_n\equiv1$ and $\beta_n\equiv0$, $\L=\L_{\rho'}$ by the uniqueness of weak limits, and thus $\seq{\QM_{T_n}^{\C_n}}$ does not converge.
\end{remark*}

\subsection{Maxima in growing dimension}
\label{sec:example_hd}

The second illustration allows the dimension to grow with $n$ and combines Theorem \ref{thm:coverage_bound} with the high-dimensional Gaussian approximation of \cite{Chernozhukov/Chetverikov/Kato:17}. For each $n\ge4$, let $X_{n,1},\ldots,X_{n,n}$ be i.i.d.\ random vectors in $\R^{p_n}$ with $p_n\ge3$. We assume that the coordinates of $X_{n,1}$ are independent with mean zero, unit variance, and sub-Gaussian norm of at most $K<\infty$. Let $m_n$ and $R_n$ be a random vector in $\R^{p_n}$ and a random $p_n\times p_n$ correlation matrix, respectively, independent of $X_{n,1},\ldots,X_{n,n}$. Take $\C_n=\sigma(m_n,R_n)$. As in Section \ref{sec:example_fd}, the observed sample consists of $Y_{n,i}$, $i=1,\ldots,n$, where
\[
	Y_{n,i}\eqdef m_n+R_n^{1/2}X_{n,i},
\]
so that $\E[Y_{n,i}\mid\C_n]=m_n$ a.s., and the parameter of interest is $\theta_n=m_n$. Similarly, $Y_n^{*}$ is drawn i.i.d.\ from the empirical distribution of $Y_n$, and $X_n^{*}$ denotes the corresponding bootstrap analog of $X_{n,1},\ldots,X_{n,n}$.

Let $Z$ be a standard normal random vector in $\R^{p_n}$, independent of everything else. We consider the statistics
\[
	T_n=\max_{1\le j\le p_n}\sqrt{n}(\widebar{Y}_{n,j}-\theta_{n,j}) \qtext{and}\quad
	T_n^{*}=\max_{1\le j\le p_n}\sqrt{n}(\widebar{Y}_{n,j}^{*}-\widebar{Y}_{n,j}),
\]
and the comparison $S_n=g(R_n,Z)$, where $g(R,x)\eqdef\max_{1\le j\le p_n}[R^{1/2}x]_j$ for a $p_n\times p_n$ correlation matrix $R$ and $x\in \R^{p_n}$. Note that
\[
	T_n=g(R_n,\sqrt{n}\widebar{X}_n) \qtext{and}\quad T_n^{*}=g(R_n,\sqrt{n}(\widebar{X}_n^{*}-\widebar{X}_n)).
\]

Fix a $p_n\times p_n$ correlation matrix $R$. The coordinates of $R^{1/2}X_{n,i}$ have the form $r^{\top}X_{n,i}$ with $\normin{r}=1$, and, by Proposition 2.7.1 in \cite{Vershynin:26:HighDim}, their sub-Gaussian norms are bounded by a constant depending only on $K$. Hence, conditions (M.1), (M.2), and (E.1) of \cite{Chernozhukov/Chetverikov/Kato:17} hold for $R^{1/2}X_{n,1},\ldots,R^{1/2}X_{n,n}$ with $b=1$ and a constant $B_n\equiv B$ depending only on $K$. Let $\XC_n\eqdef\sigma(X_{n,1},\ldots,X_{n,n})$ and $D_n\eqdef(\ln^{7}(p_nn)/n)^{1/6}$.
Since $\{x\in\R^{p_n}:\max_jx_j\le t\}$ is a hyperrectangle, Proposition 2.1 of \cite{Chernozhukov/Chetverikov/Kato:17} implies that $\dK{g(R,\sqrt{n}\widebar{X}_n),g(R,Z)\mid \H_{\varnothing}}\le C_1D_n$. Their Proposition 4.3, applied with $\alpha=n^{-1}$, yields
\[
	\DK{g(R,\sqrt{n}(\widebar{X}_n^{*}-\widebar{X}_n)),g(R,Z)\mid \XC_n,\H_{\varnothing}}\le C_2D_n
\]
with probability at least $1-n^{-1}$. The constants $C_1>0$ and $C_2>0$ depend only on $K$. Note that $(m_n,R_n)$ is independent of $X_{n,1},\ldots,X_{n,n}$ and $T_n^{*}$ depends on the data only through $Y_n$. Thus, using the identity in Remark \ref{rem:merging1}, applied to $(\omega,x)\mapsto g(R_n(\omega),x)$, we obtain
\[
	\Delta_n\le C_1D_n
	\qtext{and}\quad
	\PR{\Delta_n^{*}>C_2D_n\mid\C_n}\le n^{-1} \qtext{a.s.}
\]

Finally, $\lCZ{S_n\mid\C_n}=0$ a.s.\ because each coordinate of $R_n^{1/2}Z$ is standard normal given $\C_n$. Taking $C=C_1\vee C_2\vee 1$, we have $C\cdot D_n\ge1/n$, so that $\kappa_n\le C\cdot D_n$ a.s. Therefore, Theorem \ref{thm:coverage_bound} implies that for every $n\ge 4$,
\[
	\esssup_{\alpha\in(0,1)}\abs{\PR{T_n\le c_n^{*}(\alpha)\mid\C_n}-\alpha}\le 3C\cdot D_n \qtext{a.s.}
\]
The bound in the preceding display is nonrandom and holds for almost every realization of the location and dependence structure. It vanishes if $\ln^{7}(p_nn)=o(n)$.

\subsection{Correlations in growing dimension}
\label{sec:example_corr}

Following Remark \ref{rem:sfm_varying}, we apply the results of Section \ref{sec:linearization} to joint inference on the realized correlations in the model of Section \ref{sec:example_hd}. Suppose that $m_n=0$, $X_{n,i}\sim \nD{0}{I_{p_n}}$, and $\absin{[R_n]_{jk}}\le 1-c$ a.s.\ for all $n\ge 1$ and $1\le j<k\le p_n$ and some $c\in(0,1)$. Note that $Y_{n,1},\ldots,Y_{n,n}$ are conditionally i.i.d.\ given $\C_n$ with $Y_{n,i}\mid \C_n\sim\nD{0}{R_n}$. We assume that $\ln^{7}(p_nn)=o(n)$.

Consider $\theta_n$ and $\widehat{\theta}_n$ that consist of $[R_n]_{jk}$ and $n^{-1}\sum_{i\le n}Y_{n,ij}Y_{n,ik}$, $j\le k$, as vectors in $\R^{d_n}$ with $d_n\eqdef p_n(p_n+1)/2$, whose coordinates are indexed by such pairs $(j,k)$. Let $\widehat\theta_n^{*}$ be the bootstrap counterpart of $\widehat\theta_n$ and $\theta_n^{*}=\widehat\theta_n$. On $\Theta\eqdef\{\theta\in\R^{d_n}:\theta_{jj}>0, j=1,\ldots,p_n\}$, define $\phi_n:\Theta\to\R^{p_n(p_n-1)/2}$ by
\[
	\phi_n(\theta)_{jk}=\frac{\theta_{jk}}{\sqrt{\theta_{jj}\theta_{kk}}},\quad j<k,
\]
so that $\phi_n(\theta_n)$ and $\phi_n(\widehat{\theta}_n)$ collect the realized and sample correlations, respectively. Since $Y_{n,ij}\ne 0$ a.s., $\widehat{\theta}_n$ and $\widehat{\theta}_n^{*}$ lie in $\Theta$ a.s. We consider
\begin{align*}
	T_n&=\sqrt{n\ln{d_n}}\cdot\normin{\phi_n(\widehat{\theta}_n)-\phi_n(\theta_n)}_{\infty}
	\qtext{and} \\
	T_n^{*}&=\sqrt{n\ln{d_n}}\cdot\normin{\phi_n(\widehat{\theta}_n^{*})-\phi_n(\widehat{\theta}_n)}_{\infty},
\end{align*}
which provide simultaneous confidence intervals for the realized correlations. The factor $\sqrt{\ln{d_n}}$ does not affect the intervals or the Kolmogorov distances (Lemma \ref{lemma:scale_inv}).\footnote{
	The elements $\theta_{n,jj}$, $j=1,\ldots,p_n$ are known (equal $1$) only because the coordinates are standardized. Since the sample correlations are invariant under rescaling of the coordinates, $T_n$ and $T_n^{*}$ are unchanged if each $Y_{n,i}$ is replaced by $\Lambda_nY_{n,i}$ for an unknown $\C_n$\-/measurable diagonal matrix $\Lambda_n$ with positive diagonal entries.
}

In the notation of Section \ref{sec:linearization}, we have $\tau_n=\tau_n^{*}=\sqrt{n/\ln{d_n}}$ and $\ell$ is replaced by $\ell_n(y)=(\ln{d_n})\max_{j<k}\absin{y_{jk}}$. This normalization ensures that $\normin{\xi_n}_\infty$ is tight, since a maximum of $d_n$ centered averages is of order $\sqrt{\ln{d_n}/n}$, and it ensures that the comparison below is asymptotically equicontinuous. Each row of $D\phi_n(\theta)$ has at most three nonzero entries, which depend only on $(\theta_{jk},\theta_{jj},\theta_{kk})$. They are Lipschitz in these coordinates on $\{\theta:\normin{\theta-\theta_n}_{\infty}\le1/2\}$, with a universal Lipschitz constant. Since $\Theta$ is convex and open and the operator norm with respect to $\normin{\csdot}_{\infty}$ is the maximal absolute row sum, we have $\varpi_n(t)\le Lt$ for $t\le1/2$ and a universal constant $L$, by Remark \ref{rem:smooth}(1). Consequently, the modulus condition of Remark \ref{rem:smooth}(3) holds with $L_n=\ln{d_n}$ because $L_n\varpi_n(M/\tau_n)\le LM(\ln{d_n})^{3/2}n^{-1/2}$ for large $n$.

By \eqref{eq:linear_stat}, the linearized statistics are $T_n'=\ell_n(D\phi_n(\theta_n)\xi_n)$ and $T_n'^{*}=\ell_n(D\phi_n(\theta_n)\xi_n^{*})$. Given a $p_n\times p_n$ correlation matrix $R$, let $U_{n,i}(R)\eqdef R^{1/2}X_{n,i}$,
\begin{align*}
	\psi_{i,jk}(R)&\eqdef [U_{n,i}(R)]_j[U_{n,i}(R)]_k-[R]_{jk}\big([U_{n,i}(R)]_j^{2}+[U_{n,i}(R)]_k^{2}\big)/2, \\
	\widebar{\psi}_{n,jk}(R)&\eqdef n^{-1}\sum_{i=1}^n \psi_{i,jk}(R).
\end{align*}
We define $\psi_{i,jk}^{*}(R)$ and $\widebar{\psi}_{n,jk}^{*}(R)$ similarly with $R^{1/2}X_{n,i}^{*}$ in place of $U_{n,i}(R)$. The vectors $\psi_i(R)$ and $\widebar{\psi}_n(R)$ are formed from the elements $\psi_{i,jk}(R)$ and $\widebar{\psi}_{n,jk}(R)$, respectively, and we write $\psi_i^{*}(R)$ and $\widebar{\psi}_n^{*}(R)$ for their bootstrap counterparts.

At $\theta_n$, the nonzero entries of the $jk$-th row of $D\phi_n(\theta_n)$ are $1$, $-[R_n]_{jk}/2$, and $-[R_n]_{jk}/2$, in the coordinates $jk$, $jj$, and $kk$, and the elements of $\theta_n$ corresponding to the diagonal of $R_n$ equal $1$. Hence,
\begin{align*}
	T_n'&=\sqrt{\ln{d_n}}\cdot\norm{\sqrt{n}\widebar{\psi}_n(R_n)}_{\infty}
	\qtext{and} \\
	T_n'^{*}&=\sqrt{\ln{d_n}}\cdot\norm{\sqrt{n}(\widebar{\psi}_n^{*}(R_n)-\widebar{\psi}_n(R_n))}_{\infty}.
\end{align*}

Conditionally on $\C_n$, $\psi_{i,jk}(R_n)$ has mean zero and variance $(1-[R_n]_{jk}^2)^{2}\ge c^{2}$ a.s. Let $\Sigma_{\psi}(R)$ denote the covariance matrix of $\psi_1(R)$ and let $W(R)\eqdef{\Sigma_{\psi}}(R)^{1/2}Z$, where $Z$ is a standard normal random vector in $\R^{p_n(p_n-1)/2}$, independent of everything else. We take the comparison $S_n'=\sqrt{\ln{d_n}}\cdot \normin{W(R_n)}_{\infty}$. By Nazarov's inequality \citep[Lemma A.1]{Chernozhukov/Chetverikov/Kato:17}, applied to $(W(R_n),-W(R_n))$, we find that $\LC{\epsilon,S_n'\mid\C_n}$ is a.s.\ bounded by a constant multiple of $\epsilon$, and thus $\seq{S_n'}$ is $\CS$\-/asymptotically equicontinuous.

Next, we verify conditions (a) and (b) of Lemma \ref{lemma:bootstrap_consistency2}. Since $\absin{\psi_{i,jk}(R)}\le[R^{1/2}X_{n,i}]_j^{2}+[R^{1/2}X_{n,i}]_k^{2}$, conditions (M.1), (M.2), and (E.1) of \cite{Chernozhukov/Chetverikov/Kato:17} hold for $(\psi_i(R),-\psi_i(R))$, $1\le i\le n$, with $b=c^{2}$ and a universal constant $B_n\equiv B$. Consequently, Proposition 2.1 of \cite{Chernozhukov/Chetverikov/Kato:17} implies that, with the same $D_n$ and $\XC_n$ as in Section \ref{sec:example_hd},
\[
	\DK{\norm{\sqrt{n}\widebar{\psi}_n(R)}_{\infty},\norm{W(R)}_{\infty}\mid\H_{\varnothing}}\le C_1D_n,
\]
and Proposition 4.3 therein, applied with $\alpha=n^{-1}$, yields
\[
	\DK{\norm{\sqrt{n}(\widebar{\psi}_n^{*}(R)-\widebar\psi_n(R))}_{\infty},\norm{W(R)}_{\infty}\mid\XC_n,\H_{\varnothing}}\le C_2D_n
\]
with probability at least $1-n^{-1}$. The constants $C_1>0$ and $C_2>0$ depend only on $c$ because the dimension of $(\psi_i(R),-\psi_i(R))$ is at most $p_n^{2}$. As in Section \ref{sec:example_hd}, the identity in Remark \ref{rem:merging1}, together with Lemma \ref{lemma:scale_inv} for the factor $\sqrt{\ln{d_n}}$, yields
\begin{align*}
	&\DK{T_n',S_n'\mid\C_n}\le C_1D_n \qtext{and} \\
	&\PR{\DK{T_n'^{*},S_n'\mid\G_n,\C_n}>C_2D_n\mid\C_n}\le n^{-1} \qtext{a.s.}
\end{align*}
The same propositions, applied to the centered products $Y_{n,ij}Y_{n,ik}-[R_n]_{jk}$, yield the tightness required in Assumption \ref{assum:smooth_tightness}, because a centered normal vector $G\in\R^{d_n}$ with variances of at most $2$ satisfies $\PR{\normin{G}_{\infty}>M\sqrt{\ln{d_n}}}\le2d_n^{1-M^{2}/4}$. Consequently, the bootstrap is conditionally $\KM$\-/consistent given $\CS$ by Corollary \ref{corr:bootstrap_consistency3}.

\bibliographystyle{elsart-harv}
\bibliography{cond_bootstrap}

\end{document}